\documentclass{amsart}
\usepackage{arydshln}

\usepackage[utf8]{inputenc}
\usepackage[T1]{fontenc}

\usepackage{amsmath,amsthm,amsopn,amstext,amscd,amsfonts,amssymb,mathrsfs,mathtools}
\usepackage{dsfont}
\usepackage{comment}
\usepackage{xcolor}
\usepackage{float}
\usepackage[active]{srcltx}
\usepackage{graphicx, mathdots}

  \DeclareMathOperator{\sgn}{sgn}
  \usepackage{mathtools}

\newtheorem{theorem}{\sc Theorem}[section]

\newtheorem{lemma}{\sc Lemma }[section]
\newtheorem{eje}{\sc Example }[section]

\newtheorem{prop}{\sc Proposition}[section]
\newtheorem{remark}{\sc Remark}[section]
\usepackage[active]{srcltx}

\newcommand{\dps}{\displaystyle}
\newcommand{\supp}{\operatorname{supp}}
\usepackage{graphicx, epsfig, subfig}
\usepackage{lscape}
\usepackage{rotating}
\usepackage{caption}
\usepackage{multicol}
\usepackage{colortbl}
\usepackage{nicematrix}
\usepackage{hyperref}

\begin{document}
\title[The multiple Markov theorem on Angelesco sets]{The multiple Markov theorem on Angelesco sets}
\author{K. Castillo}
\address{CMUC, Department of Mathematics, University of Coimbra, 3000-143 Coimbra, Portugal}
\email{kenier@mat.uc.pt}
\author{G. Gordillo-N\'u\~nez}
\address{CMUC, Department of Mathematics, University of Coimbra, 3000-143 Coimbra, Portugal}
\email{up202310693@up.pt}
\subjclass[2010]{26C10, 33C47}
\date{\today}
\keywords{}
\begin{abstract}
Addressing a problem described by Ismail as a ``worthwhile research
project'', this note extends Markov's theorem to multiple orthogonal
polynomials on Angelesco sets. The interaction $\mathcal{Z}$-matrix is shown
to be a nonsingular $\mathcal{M}$-matrix; hence the theorem
holds for positive Borel measures satisfying the natural support condition,
including finite supports, without restriction on
the number of intervals. Beyond separated supports, a dual variation identity
gives necessary and sufficient tail conditions. An explicit AT example
disproves the componentwise extension, while for Nikishin systems we prove
positive results for certain deformations of the first generator and give a
counterexample for an internal generator.
\end{abstract}
\maketitle
\section{Introduction\label{intro}}
In the third volume of {\em Mathematics of the 19th Century}, edited by Kolmogorov and Yushkevich (see \cite[p. 52]{KY98}\footnote{Originally published as “Matematika XIX veka: Chebyshevskoe napravlenie v teorii funktsii. Obyknovennye differentsial'nye uravneniia. Variatsionnoe ischislenie. Teoriia konechnykh raznostei,” Izdatel’stvo “Nauka”, Moskva, 1987.}), three papers by A. Markov addressing the behaviour of the zeros of the non-constant elements in a sequence of orthogonal polynomials on the real line (OPRL), $(p_n)_{n=0}^\infty$, with respect to a parameter, $t$, varying within a real open interval, are discussed. Among these theorems, one stands out as the most elegant of Markov’s results on orthogonal polynomials, highlighting the pivotal role of the orthogonality weight function, $\omega$, defined on a (not necessarily finite) real interval $(a, b)$ (see \cite[p. 178]{M86}):

\vspace{2mm}
\begin{changemargin}{0.8cm}{0.8cm}
\emph{Under suitable conditions, the zeros of $p_n(x; t)$ are increasing (respectively, decreasing) functions of $t$, provided that
$$
\frac{1}{\omega(x; t)}\frac{\partial \omega}{\partial t}(x; t)
$$
is an increasing (respectively, decreasing)\footnote{When strictly monotonic behaviour is present in the hypotheses, the zeros also become strictly monotonic functions with respect to the parameter.} function of $x$ on $(a, b)$.}
\end{changemargin}
\vspace{2mm}
From Szeg\H{o}'s seminal 1939 book (see \cite[Theorem 6.12.1]{S75}) to Ismail’s monograph (see \cite[Theorem 7.1.1]{I05}), and through Freud’s book (see \cite[Problems 14, 15, and 16, pp. 133–134]{G71}), Markov’s theorem has frequently appeared in the literature. Notably, Askey and Wilson used a direct consequence of Markov’s theorem, originally traced back to Szeg\H{o} (see \cite[Theorem 6.12.2]{S75}), to explore the monotonicity of the zeros of what are now known as the Askey-Wilson polynomials (see \cite[Section 7]{AW85}). 

The monotonic behaviour of the zeros of classical OPRL with respect to their parameters can often be deduced from their specific properties, including Rodrigues’ formula, differential or difference equations, hypergeometric character, among others. In many instances, the simplicity of the problem is such that the expected result can almost be derived directly. For example, in the case of Laguerre polynomials, the problem can be reformulated as an eigenvalue problem, with the result following from the application of Hadamard's first variation formula (see, for instance,  \cite[(1.73)]{Tao} or \cite[Proposition 3.1]{CZ22}). However, it is its generality and simplicity that, as Szeg\H{o} described, makes Markov’s theorem ``an important statement'', or, in the words of Ismail, ``an extremely useful theorem''. The Laguerre polynomial of degree $n$, $L^{(\alpha)}_n(x)$, with $\alpha \in (-1, \infty)$, satisfies the following orthogonality conditions:
$$
\int_{0}^\infty  x^{j} L^{(\alpha)}_n(x)e^{-x} x^\alpha \,\mathrm{d}x = 0,\quad j = 0, 1, \dots, n-1.
$$
As previously mentioned, it would require little effort to observe that their zeros, whose weight function is given by  
$$
\omega(x; \alpha) = e^{-x} x^{\alpha},
$$
are increasing functions of $\alpha$ on $(-1,\infty)$. However, according to Markov's theorem, this property effortlessly emerges because
\begin{align}\label{LaguerreMarkov}
\frac{1}{\omega(x; \alpha)}
\frac{\partial \omega}{\partial \alpha}(x; \alpha)=\log x
\end{align}
is an increasing function of $x$ on $(0, \infty)$. Recently, relying solely on orthogonality conditions, it was proved in \cite{C22} that Markov’s theorem remains essentially valid for weight functions defined on the unit circle, opening promising avenues for extending Markov’s theorem to other contexts, such as multiple orthogonal polynomials, where applications demand such results.

Over the past four decades, multiple orthogonal polynomials have garnered significant interest due to their extensive theoretical development and relevance in various mathematical contexts. The introduction of \cite{D23} provides an interesting and up-to-date discussion in this regard. However, a generalisation of Markov’s theorem to the setting of “several weight functions” has yet to be established, complicating the study of the monotonicity of their zeros. Although it had long been recognised as a significant challenge, in 2005, Ismail explicitly included this question in his collection of open problems and conjectures, which encompasses some of the most challenging problems in the field of Orthogonal Polynomials and Special Functions (see \cite[Problem 24.1.5]{I05}):
\vspace{2mm}
\begin{changemargin}{0.8cm}{0.8cm}
\emph{“There is no study of zeros of general or special systems of multiple orthogonal polynomials available. An extension of Theorem 7.1.1 $[$Markov’s theorem$]$ to multiple orthogonal polynomials would be a worthwhile research project.”}
\end{changemargin}
\vspace{2mm}
We prove that Markov's theorem survives on Angelesco sets. More precisely,
if the logarithmic derivative of each weight with respect to the parameter
is non-decreasing on the interval carrying its measure, then every zero
of the multiple orthogonal polynomial is non-decreasing in that parameter.
The analogous statement holds with non-increasing in place of non-decreasing.
This gives a general answer to Ismail's
extension problem for Angelesco systems. Beyond separated supports, a dual
variation identity yields necessary and sufficient tail conditions for the
independent and common deformation classes. The AT counterexample shows
that the componentwise condition alone is insufficient in general, while
the Nikishin results establish further positive cases and an obstruction
for deformations of an internal generator.

\begin{eje}\label{Example1}
Let $a>0$ and $b>0$, and let $\{\nu_1, \nu_2\}$ denote a set of two
nontrivial positive Borel measures on $\mathbb{R}$ defined as
$$
\mathrm{d}\nu_k(x; a_k, b_k, c_k)
=\omega(x; a_k, b_k, c_k)\,\mathrm{d}x,\quad k = 1, 2.
$$
We suppress the endpoints $a$ and $b$ from the parameter lists; when either
is varied, all other parameters are fixed.
The weight functions are expressed as
$$
\omega(x; a_k, b_k, c_k) = (x+a)^{a_k} |x|^{b_k} (b-x)^{c_k},\quad a_k, b_k, c_k \in (-1, \infty).
$$
Assume that $\nu_1$ is concentrated on the interval $(-a,0)$, while
$\nu_2$ is concentrated on $(0,b)$, as illustrated in the following figure:
\begin{figure}[H]
\centering
\includegraphics[width=8cm]{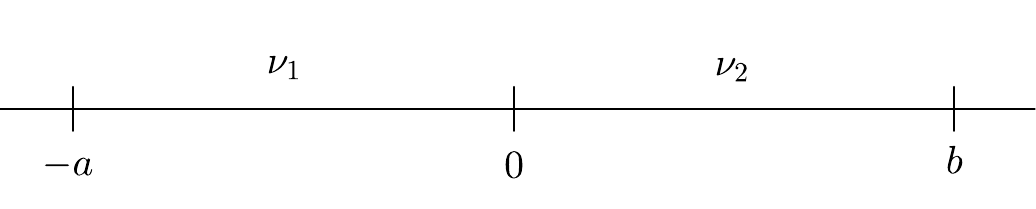}
\caption{The intervals carrying $\nu_1$ and $\nu_2$.}\label{fig:two-supports}
\end{figure} 
\noindent (After the changes $u=(x+a)/a$ on $(-a,0)$ and $u=x/b$ on
$(0,b)$, these are, up to positive constants, Magnus's generalised Jacobi
weights with three factors, related to the sixth Painlev\'e equation; see
\cite[Example~1]{M95}.)

Let positive integers $n_1,n_2$ be fixed and put $n=n_1+n_2$. Since the
measures $\nu_1$ and $\nu_2$ are concentrated on pairwise disjoint open real intervals,
there exists a MOPRL of degree $n$,
$P_n^{(a_1,b_1,c_1,a_2,b_2,c_2)}$, which satisfies the following
orthogonality conditions, with $0\le j\le n_k-1$ in the condition
corresponding to $\nu_k$:
\begin{align*}
\int_{-a}^0 x^j P_n^{(a_1, b_1, c_1, a_2, b_2, c_2)}(x) \, \omega(x; a_1, b_1, c_1) \,\mathrm{d}x &= 0, \\[7pt]
\int_0^b x^j P_n^{(a_1, b_1, c_1, a_2, b_2, c_2)}(x) \, \omega(x; a_2, b_2, c_2) \,\mathrm{d}x &= 0.
\end{align*}
A straightforward calculation shows that  
\begin{align*}
\frac{1}{\omega(x; a_k, b_k, c_k)}
\frac{\partial \omega}{\partial a_k}(x; a_k, b_k, c_k)
&=\log(x+a),\\[7pt]
\frac{1}{\omega(x; a_k, b_k, c_k)}
\frac{\partial \omega}{\partial b_k}(x; a_k, b_k, c_k)
&=\log |x|,\\[7pt]
\frac{1}{\omega(x; a_k, b_k, c_k)}
\frac{\partial \omega}{\partial c_k}(x; a_k, b_k, c_k)
&=\log(b-x),\quad k=1,2.
\end{align*}
The first and third functions are, respectively, increasing and decreasing on
$(-a,b)$, while $\log|x|$ decreases on $(-a,0)$ and increases on $(0,b)$.
The weight indexed by $k$ is independent of the parameters carrying the other
index. It will follow that the
zeros of $P_n^{(a_1,b_1,c_1,a_2,b_2,c_2)}$ are increasing functions of
$a_1,a_2,b_2$ and decreasing functions of $b_1,c_1,c_2$ throughout the natural
range $a_k,b_k,c_k>-1$ (see Theorem~\ref{maintheorem} below).  In the particular case $a_1=a_2$, $b_1=b_2$,
$c_1=c_2$, and $b=1$, one obtains the Jacobi--Angelesco polynomials (see
\cite[Section 23.3.1]{I05}).  If in addition $n_1=n_2$, their Rodrigues formula
allows more elementary methods; see \cite{D17,MM24} and \cite{L22}.
\end{eje}

Since it is not necessary to know the explicit form of MOPRL on Angelesco sets, nor any of their specific properties, to predict the possible monotonic behaviour of their zeros, we have complete freedom in choosing examples. In this sense, and as an illustration of the strength of the main theorem to be proved in this work, a final simple example is presented, whose result is far from elementary to predict, even computationally.

\begin{eje}
    Let $\{\nu_1,\nu_2,\nu_3\}$ denote a set of three nontrivial positive Borel measures on $\mathbb{R}$ defined as   
$$
    \mathrm{d}\nu_k(x; a_k, b_k)=\omega(x; a_k, b_k)\,\mathrm{d}x,\quad k = 1, 2, 3,
$$
with the weight functions expressed as
    \begin{align*}
        \omega(x; a_1, b_1) &= \left(-1-x\right)^{a_1}e^{b_1\, x},\quad a_1 \in (-1, \infty),\quad b_1\in(0, \infty),\\[7pt]
        \omega(x; a_2, b_2) &= (1+x)^{a_2}(-x)^{b_2},\quad a_2, b_2 \in (-1, \infty), \\[7pt]
        \omega(x; a_3, b_3) &= e^{-a_3\,x^3+b_3x},\quad
        a_3 \in (0, \infty),\quad b_3\in\mathbb{R}.
    \end{align*}
Assume that $\nu_1$, $\nu_2$, and $\nu_3$ are concentrated, respectively,
on $(-\infty,-1)$, $(-1,0)$, and $(0,\infty)$, as illustrated in the
following figure:
       \begin{figure}[H]
\centering
\includegraphics[width=8cm]{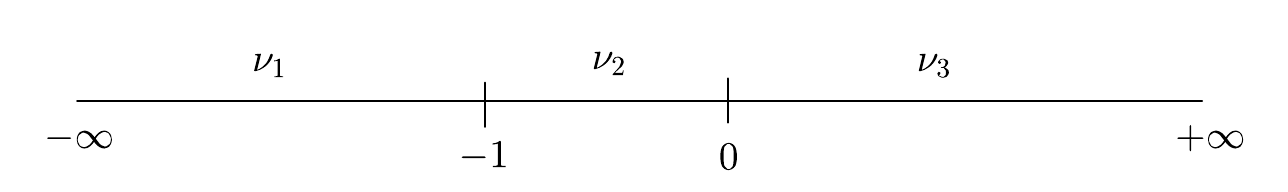}
\caption{The intervals carrying $\nu_1$, $\nu_2$, and $\nu_3$.}\label{fig:three-supports}
\end{figure}
\noindent (The first two measures are classical Laguerre and Jacobi weights.
The third is a cubic exponential weight on the positive half-line.)

Let positive integers $n_1,n_2,n_3$ be fixed and put
$n=n_1+n_2+n_3$. Since the measures $\nu_1$, $\nu_2$, and $\nu_3$ are
concentrated on pairwise disjoint open real intervals, there exists a MOPRL of degree $n$,
$P_n^{(a_1,b_1,a_2,b_2,a_3,b_3)}$, which satisfies the following
orthogonality conditions, with $0\le j\le n_k-1$ in the condition
corresponding to $\nu_k$:
\begin{align*}
\int_{-\infty}^{-1} x^j P_n^{(a_1, b_1, a_2, b_2, a_3, b_3)}(x) \, \omega(x; a_1, b_1) \,\mathrm{d}x &= 0, \\[7pt]
\int_{-1}^0 x^j P_n^{(a_1, b_1, a_2, b_2, a_3, b_3)}(x) \, \omega(x; a_2, b_2) \,\mathrm{d}x &= 0, \\[7pt]
\int_0^\infty x^j P_n^{(a_1, b_1, a_2, b_2, a_3, b_3)}(x) \, \omega(x; a_3, b_3) \,\mathrm{d}x &= 0.
\end{align*}
 A straightforward calculation shows that  
\begin{align*}
\frac{1}{\omega(x; a_1, b_1)}
\frac{\partial \omega}{\partial a_1}(x; a_1, b_1)
&=\log |1+x|,\\[7pt]
\frac{1}{\omega(x; a_2, b_2)}
\frac{\partial \omega}{\partial a_2}(x; a_2, b_2)
&=\log |1+x|,\\[7pt]
\frac{1}{\omega(x; a_1, b_1)}
\frac{\partial \omega}{\partial b_1}(x; a_1, b_1)
&=x,\\[7pt]
\frac{1}{\omega(x; a_3, b_3)}
\frac{\partial \omega}{\partial b_3}(x; a_3, b_3)
&=x,\\[7pt]
\frac{1}{\omega(x; a_2, b_2)}
\frac{\partial \omega}{\partial b_2}(x; a_2, b_2)
&=\log |x|,\\[7pt]
\frac{1}{\omega(x; a_3, b_3)}
\frac{\partial \omega}{\partial a_3}(x; a_3, b_3)
&=-x^3.
\end{align*}
Thus the corresponding logarithmic derivatives for $a_1,b_2,a_3$ are decreasing, while those for
$b_1,a_2,b_3$ are increasing, on their respective supports.  We shall conclude
that the zeros of $P_n^{(a_1,b_1,a_2,b_2,a_3,b_3)}$ increase with
$b_1,a_2,b_3$ and decrease with $a_1,b_2,a_3$ throughout the parameter ranges
stated above  (see Theorem~\ref{maintheorem} below).
\end{eje}
 
We next establish differentiability of the zeros and prove the multiple
Markov theorem and its extensions.
\section{Preliminary results}
Let $\{\nu_1,\dots,\nu_m\}$ be a nonempty set of nontrivial positive Borel
measures on $\mathbb{R}$, each having finite moments. Fix a multi-index
$(n_1,\dots,n_m)$ with positive entries and put
$n=n_1+\cdots+n_m$. A component with $n_k=0$ imposes no orthogonality
condition and may simply be omitted, so this convention entails no loss.
For $r=0,1,\ldots$, we write $\mathbb{P}_r$ for the space of real polynomials of
degree at most $r$.
A monic polynomial $P_n$ of degree $n$ associated with
the measures and this multi-index is called a MOPRL if it satisfies the
following orthogonality conditions (see \cite[Definition 1.1]{A98}):
\begin{align}\label{orthomultiple}
\int x^j P_n(x)\,\mathrm{d}\nu_k(x) = 0,\quad
j = 0, 1, \dots, n_k - 1,\quad k = 1, 2, \dots, m.
\end{align}
In some contexts, these polynomials are referred to as type II MOPRL. Their
existence is not guaranteed for an arbitrary system and an arbitrary
multi-index. Under the standard infinite-support Angelesco hypotheses,
existence and uniqueness are classical (see \cite[p. 135]{NS91} and
\cite[Corollary 1.5]{A98}). We use the same terminology, including for
finite supports, when there are pairwise disjoint open real intervals
$\Delta_k$ on which the respective measures are concentrated and
$\#\supp\nu_k\ge n_k$. For
completeness, suppose that a nonzero polynomial $S$ of degree at most $n-1$
satisfies the homogeneous conditions. In $\Delta_k$ it either vanishes
at every support point of the measure lying there or has at least $n_k$
sign changes. Otherwise,
multiplication by the product of its sign-change points, with a suitable
choice of sign, produces a polynomial of degree below $n_k$ whose product
with $S$ is non-negative and nonzero on a set of positive measure, contrary
to orthogonality. Since the measure is concentrated on $\Delta_k$ and
$\#\supp\nu_k\ge n_k$, either alternative gives at least $n_k$ distinct
zeros in $\Delta_k$. The disjoint intervals contradict
$\deg S\le n-1$. Thus the moment matrix is nonsingular, and the same argument
applied to $P_n$ gives exactly $n_k$ simple zeros in each block. In this case,
the set $\{\nu_1,\dots,\nu_m\}$ is referred to as an
Angelesco set or Angelesco system (see \cite[Definition 1.3]{A98}). Henceforth,
in the Angelesco arguments of the main section, $P_n$ refers to the polynomial
associated with the displayed measures and the fixed multi-index. The notation is particularly
well-suited to the problem at hand and remains clear for a reader unfamiliar
with MOPRL.

\begin{prop}\label{diff}
Let $I$ be an open interval containing $t_0$, let
$n_1,\ldots,n_m$ be positive integers, and put $n=n_1+\cdots+n_m$.
For $k=1,\ldots,m$, let $\mu_k$ be a fixed positive Borel measure and set
$$
\mathrm{d}\nu_k(x; t)=\omega_k(x; t)\,\mathrm{d}\mu_k(x),
$$
where $\omega_k(x; t)>0$ for $\mu_k$-almost every $x$ and
$\mu_k(\mathbb{R}\setminus(a_k,b_k))=0$. For every $k=1,\ldots,m$ and
$0\le r\le n+n_k-1$, assume that
$$
c_r^{(k)}(t)=\int x^r\omega_k(x; t)\,\mathrm{d}\mu_k(x)
$$
is $C^1$ on $I$. Suppose also that, for every $k=1,\ldots,m$ and $t\in I$,
there is a Borel function $x\longmapsto\phi_k(x; t)$ such that
\begin{align}\label{eq:representing-moments}
\int(1+|x|^{n+n_k-1})|\phi_k(x; t)|
\,\mathrm{d}\nu_k(x; t)&<\infty,\\[7pt]
\frac{\mathrm{d}c_r^{(k)}}{\mathrm{d}t}(t)
 &=\int x^r\phi_k(x; t)\,\mathrm{d}\nu_k(x; t).
\end{align}
We shall call any such function an admissible representing function. If
$\frac{\partial\omega_k}{\partial t}(x; t)$ exists, differentiation under
the integral sign is justified, and the displayed integrability condition
holds, then the logarithmic derivative
$$
 \frac{1}{\omega_k(x; t)}
 \frac{\partial\omega_k}{\partial t}(x; t)
$$
is an admissible representing function.
For every parameter value at which the type II problem for the measures
$\nu_1,\ldots,\nu_m$ and the multi-index $(n_1,\ldots,n_m)$ is normal,
denote its monic polynomial by $P_n(x; t)$.
If the multiple
orthogonality problem is normal near $t_0$ and $y_0$ is a simple zero of
$P_n(x; t_0)$, then there exist
$\epsilon,\delta>0$ and a unique $C^1$ function
$$
y:(t_0-\epsilon,t_0+\epsilon)\longrightarrow(y_0-\delta,y_0+\delta)
$$
such that $y(t_0)=y_0$ and $P_n(y(t); t)=0$.
\end{prop}
\begin{proof}
Write $P_n(x; t)=x^n+\sum_{r=0}^{n-1}p_r(t)x^r$ and set
$\mathbf{p}(t)=(p_0(t),\ldots,p_{n-1}(t))^{\mathrm T}$. The orthogonality conditions
form the linear system
\begin{align*}
\mathbf{M}(t)\mathbf{p}(t)&=-\mathbf{h}(t),\\[7pt]
m_{(k,j),r}(t)&=c_{j+r}^{(k)}(t),\\[7pt]
h_{(k,j)}(t)&=c_{j+n}^{(k)}(t),
\end{align*}
where $m_{(k,j),r}(t)$ is the entry of $\mathbf{M}(t)$ in row $(k,j)$ and
column $r$, and $h_{(k,j)}(t)$ is the corresponding component of
$\mathbf{h}(t)$; here $1\le k\le m$, $0\le j<n_k$, and $0\le r<n$.
Normality is precisely
$\det\mathbf{M}(t)\ne0$.  Thus Cramer's rule and the hypotheses on the moments show
that all coefficients $p_r$ are $C^1$.  Since
$\frac{\partial P_n}{\partial x}(y_0; t_0)\ne0$, the $C^1$ implicit
function theorem gives the
assertion.
\end{proof}

The proposition does not require separated supports; separation enters only
in the sign analysis. Section~\ref{sec:dual} gives the general identity.

\section{Main results}
The proof of the main result will follow directly from a series of lemmas. We highlight that many of these lemmas, and consequently the monotonicity of the zeros of MOPRL on Angelesco sets, emerge exclusively from the orthogonality properties \eqref{orthomultiple}. 
After relabelling, choose pairwise disjoint open intervals
$\Delta_k=(a_k,b_k)$ on which the measures $\mu_k$ are concentrated and
ordered so that
$b_k\le a_{k+1}$.

Throughout, matrices are denoted by bold upper-case letters and vectors by
bold lower-case letters; their entries and components are written in ordinary
type and carry indices. Inequalities between matrices are understood entrywise
and inequalities between vectors componentwise. Thus $\mathbf{C}\ge0$ means that
every entry of the matrix $\mathbf{C}$ is non-negative, while $\mathbf{C}>0$
means that every entry is positive. For a square matrix $\mathbf{C}$,
$\rho(\mathbf{C})$ denotes its spectral radius.
Recall that a $\mathcal{Z}$-matrix is a real square
matrix with non-positive off-diagonal entries. A nonsingular
$\mathcal{M}$-matrix is a $\mathcal{Z}$-matrix that is invertible and whose
inverse is non-negative \cite[Section~2.3 and Theorem~5.2.1]{JST20}.

\begin{lemma}\label{lemma1}
Assume the hypotheses and notation of Proposition~\ref{diff}. Suppose that
the intervals $(a_k,b_k)$ are pairwise disjoint and that
$\#\supp\mu_k\ge n_k$ for every $k$. Thus the measures form an Angelesco
system at every parameter value. Let $Q_{n_k}(x; t)$ denote the monic
polynomial whose zeros, $x_{k,1}(t)<\cdots<x_{k,n_k}(t)$, are precisely the
zeros of $P_n(x; t)$ in $(a_k,b_k)$. The notation $Q_{n_k}$ retains the
block label $k$, even when different blocks have the same degree.
Define the block matrix
$$
\mathbf{A}(t) =
\begin{pmatrix}
\mathbf{I}_{n_1} & \mathbf{A}_{1,2}(t) & \cdots & \mathbf{A}_{1,m-1}(t) & \mathbf{A}_{1,m}(t) \\[7pt]
\mathbf{A}_{2,1}(t) & \mathbf{I}_{n_2} & \cdots & \mathbf{A}_{2,m-1}(t) & \mathbf{A}_{2,m}(t) \\[7pt]
\vdots & \vdots & \ddots & \vdots & \vdots \\[7pt]
\mathbf{A}_{m-1,1}(t) & \mathbf{A}_{m-1,2}(t) & \cdots & \mathbf{I}_{n_{m-1}} & \mathbf{A}_{m-1,m}(t) \\[7pt]
\mathbf{A}_{m,1}(t) & \mathbf{A}_{m,2}(t) & \cdots & \mathbf{A}_{m,m-1}(t) & \mathbf{I}_{n_m}
\end{pmatrix},
$$
where $\mathbf{I}_r$ denotes the identity matrix of order $r$, and
$\mathbf{A}_{k,l}(t)$ is defined as
$$
\mathbf{A}_{k,l}(t) =
\begin{pmatrix}
a^{(k,l)}_{1,1}(t) & a^{(k,l)}_{1,2}(t) & \cdots & a^{(k,l)}_{1,n_l-1}(t) & a^{(k,l)}_{1,n_l}(t) \\[7pt]
a^{(k,l)}_{2,1}(t) & a^{(k,l)}_{2,2}(t) & \cdots & a^{(k,l)}_{2,n_l-1}(t) & a^{(k,l)}_{2,n_l}(t) \\[7pt]
\vdots & \vdots & \ddots & \vdots & \vdots \\[7pt]
a^{(k,l)}_{n_k-1,1}(t) & a^{(k,l)}_{n_k-1,2}(t) & \cdots & a^{(k,l)}_{n_k-1,n_l-1}(t) & a^{(k,l)}_{n_k-1,n_l}(t) \\[7pt]
a^{(k,l)}_{n_k,1}(t) & a^{(k,l)}_{n_k,2}(t) & \cdots & a^{(k,l)}_{n_k,n_l-1}(t) & a^{(k,l)}_{n_k,n_l}(t)
\end{pmatrix}.
$$
For $k,l=1,\dots,m$ with $k\ne l$, $i=1,\dots,n_k$, and
$j=1,\dots,n_l$, put
$$
d_i^{(k)}(t)=\int
\frac{Q_{n_k}(x; t)}{x-x_{k,i}(t)}
\frac{P_n(x; t)}{x-x_{k,i}(t)}
\,\mathrm{d}\nu_k(x; t),
$$
and define
\begin{align*}
a^{(k,l)}_{i,j}(t)
&=\frac{1}{d_i^{(k)}(t)}\int
\frac{Q_{n_k}(x; t)}{x-x_{k,i}(t)}
\frac{P_n(x; t)}{x-x_{l,j}(t)}
\,\mathrm{d}\nu_k(x; t).
\end{align*}
Every quotient of a polynomial by one of its displayed zeros is understood
as the polynomial obtained after exact division, also at that zero.
Then there exists a neighbourhood of $t_0$ in which $\mathbf{A}(t)$ is a
$\mathcal{Z}$-matrix.
\end{lemma}
\begin{proof}
Since every zero of $Q_{n_l}(x; t)$ belongs to $\Delta_l$, the polynomial
$$
R_k(x; t)=\prod_{\substack{l=1 \\[7pt] l \neq k}}^{m} Q_{n_l}(x; t),
$$
has no zero and has constant sign on $\Delta_k$. Moreover,
$$
d_i^{(k)}(t)=\int R_k(x; t)
\left(\frac{Q_{n_k}(x; t)}{x-x_{k,i}(t)}\right)^2
\,\mathrm{d}\nu_k(x; t).
$$
The polynomial $Q_{n_k}(x; t)/(x-x_{k,i}(t))$ has only $n_k-1$
zeros, whereas the support of the $k$th measure at $t$ contains at least
$n_k$ points. Hence the last integral is nonzero and
$$
\sgn d_{i}^{(k)}(t) = \sgn\, R_k(x; t)
$$
on $\Delta_k$. Put $f(x; t)=1/(x-x_{l,j}(t))$. Orthogonality gives
$$
\frac{1}{x_{k,i}(t) - x_{l,j}(t)} \int \frac{Q_{n_k}(x; t)P_n(x; t)}{x - x_{k,i}(t)}\,\mathrm{d}\nu_k(x; t) = 0.
$$
After subtracting this zero integral from the numerator that defines
$a^{(k,l)}_{i,j}(t)d_i^{(k)}(t)$, we obtain
\begin{align*}
a^{(k,l)}_{i,j}(t)d_i^{(k)}(t)
&=\int R_k(x; t)
\left(\frac{Q_{n_k}(x; t)}{x-x_{k,i}(t)}\right)^2
\\[7pt]
&\hspace{12mm}{}\times(x-x_{k,i}(t))
\bigl(f(x; t)-f(x_{k,i}(t); t)\bigr)
\,\mathrm{d}\nu_k(x; t).
\end{align*}
Here the polynomial quotients are understood after cancellation, and
$$
(x-x_{k,i}(t))\bigl(f(x; t)-f(x_{k,i}(t); t)\bigr)
=-\frac{(x-x_{k,i}(t))^2}
{(x-x_{l,j}(t))(x_{k,i}(t)-x_{l,j}(t))}\le0
$$
on $\Delta_k$. The formula is also valid at $x=x_{k,i}(t)$, where both
sides are zero. The two factors in the denominator have the same sign
because $\Delta_k$ and $\Delta_l$ are disjoint. Therefore every
off-diagonal entry satisfies $a^{(k,l)}_{i,j}(t)\le0$. If the support of
the $k$th measure at $t$ contains more than $n_k$ points, it is not
contained in the zeros of $Q_{n_k}(x; t)$, and the inequality is strict.
\end{proof}

\begin{lemma}\label{Ax=b}
Assume the hypotheses and notation of Lemma \ref{lemma1}. Define the vectors
\begin{align*}
\mathbf{x}(t)&=\bigl(x_{1,1}(t),\ldots,x_{1,n_1}(t),\ldots,
x_{m,1}(t),\ldots,x_{m,n_m}(t)\bigr)^{\mathrm T},\\[7pt]
\mathbf{b}(t)&=\bigl(b^{(1)}_1(t),\ldots,b^{(1)}_{n_1}(t),\ldots,
b^{(m)}_1(t),\ldots,b^{(m)}_{n_m}(t)\bigr)^{\mathrm T},
\end{align*}
where, for $k=1,\dots,m$ and $i=1,\dots,n_k$,
$$
b^{(k)}_i(t)=\frac{1}{d_i^{(k)}(t)}\int
Q_{n_k}(x; t)\frac{P_n(x; t)}{x-x_{k,i}(t)}
\phi_k(x; t)\,\mathrm{d}\nu_k(x; t).
$$
In this block order, the component with index $(k,i)$ occupies position
$n_1+\cdots+n_{k-1}+i$, with the preceding sum understood as zero when $k=1$.
Then there exists a neighbourhood of $t_0$ in which
$\mathbf{A}(t)\frac{\mathrm{d}\mathbf{x}}{\mathrm{d}t}(t)=\mathbf{b}(t)$.
\end{lemma}

\begin{proof}
Differentiate first the monomial orthogonality relations. At the fixed value
of $t$, freeze the coefficients of the polynomial
$h(x; t)=Q_{n_k}(x; t)/(x-x_{k,i}(t))$ and take the corresponding linear
combination of the differentiated relations. Thus no derivative of
$h(x; t)$ is introduced, and
\begin{align*}
0={}&\int h(x; t)\frac{\partial Q_{n_k}}{\partial t}(x; t)
R_k(x; t)\,\mathrm{d}\nu_k(x; t)\\[7pt]
&+\int h(x; t)Q_{n_k}(x; t)
\frac{\partial R_k}{\partial t}(x; t)\,\mathrm{d}\nu_k(x; t)\\[7pt]
&+\int h(x; t)Q_{n_k}(x; t)R_k(x; t)
\phi_k(x; t)\,\mathrm{d}\nu_k(x; t).
\end{align*}
For every polynomial $g$ of degree at most $n_k-1$, orthogonality also gives
$$
\int\frac{P_n(x; t)}{x-x_{k,i}(t)}g(x)\,\mathrm{d}\nu_k(x; t)
 =g(x_{k,i}(t))\int\frac{P_n(x; t)}{x-x_{k,i}(t)}
 \,\mathrm{d}\nu_k(x; t).
$$
Indeed, after bringing the right-hand side to the left, the remaining
multiplier of $P_n(x; t)$ is
$\bigl(g(x)-g(x_{k,i}(t))\bigr)/(x-x_{k,i}(t))$, a polynomial of degree at
most $n_k-2$; it is the zero polynomial when $n_k=1$. Apply the identity
first to $g(x)=\frac{\partial Q_{n_k}}{\partial t}(x; t)$ and then to
$g(x)=Q_{n_k}(x; t)/(x-x_{k,i}(t))$. Since
$\frac{\partial Q_{n_k}}{\partial t}(x_{k,i}(t); t)
=-\frac{\mathrm{d}x_{k,i}}{\mathrm{d}t}(t)
\frac{\partial Q_{n_k}}{\partial x}(x_{k,i}(t); t)$ and
$$
 \frac{1}{R_k(x; t)}\frac{\partial R_k}{\partial t}(x; t)
 =-\sum_{l\ne k}\sum_{j=1}^{n_l}
 \frac{1}{x-x_{l,j}(t)}\frac{\mathrm{d}x_{l,j}}{\mathrm{d}t}(t),
$$
the three integrals above are, respectively,
\begin{align*}
-\frac{\mathrm{d}x_{k,i}}{\mathrm{d}t}(t)d_i^{(k)}(t),\quad
-d_i^{(k)}(t)\sum_{l\ne k}\sum_{j=1}^{n_l}
a_{i,j}^{(k,l)}(t)\frac{\mathrm{d}x_{l,j}}{\mathrm{d}t}(t),\quad
d_i^{(k)}(t)b_i^{(k)}(t).
\end{align*}
The second equality follows by substituting the displayed formula for
$\frac{\partial R_k}{\partial t}(x; t)$ and the definition of
$a_{i,j}^{(k,l)}(t)$; the third is the definition of $b_i^{(k)}(t)$.
For the first equality we have also used
$$
d_i^{(k)}(t)=
\frac{\partial Q_{n_k}}{\partial x}(x_{k,i}(t); t)
\int\frac{P_n(x; t)}{x-x_{k,i}(t)}\,\mathrm{d}\nu_k(x; t),
$$
which is the preceding identity with
$g(x)=Q_{n_k}(x; t)/(x-x_{k,i}(t))$. Their sum is zero, and division by the
nonzero number $d_i^{(k)}(t)$ gives
$$
 \frac{\mathrm{d}x_{k,i}}{\mathrm{d}t}(t)
 +\sum_{l\ne k}\sum_{j=1}^{n_l}a_{i,j}^{(k,l)}(t)
 \frac{\mathrm{d}x_{l,j}}{\mathrm{d}t}(t)=b_i^{(k)}(t).
$$
These are precisely the rows of
$\mathbf{A}(t)\frac{\mathrm{d}\mathbf{x}}{\mathrm{d}t}(t)=\mathbf{b}(t)$.
\end{proof}

The following lemma contains Markov's sign argument.

\begin{lemma}\label{aux1}
Assume the hypotheses and notation of Lemma \ref{Ax=b}. Then $b^{(k)}_i(t)$
is non-negative in a neighbourhood of $t_0$, provided that the functions
$x\longmapsto\phi_k(x; t)$ have non-decreasing representatives on
$\Delta_k$.
\end{lemma}
\begin{proof}
Choose a non-decreasing representative of $x\longmapsto\phi_k(x; t)$ on
$\Delta_k$.
With $c_{k,i}(t)=\phi_k(x_{k,i}(t); t)$, orthogonality permits us to subtract the
constant $c_{k,i}(t)$ and gives
\begin{align*}
b_i^{(k)}(t)
 &=\frac{1}{d_i^{(k)}(t)}\int R_k(x; t)
 \left(\frac{Q_{n_k}(x; t)}{x-x_{k,i}(t)}\right)^2
 \\[7pt]
 &\hspace{12mm}{}\times(x-x_{k,i}(t))
 \bigl(\phi_k(x; t)-c_{k,i}(t)\bigr)
 \,\mathrm{d}\nu_k(x; t).
\end{align*}
The last product is defined to be zero at $x=x_{k,i}(t)$. It is
non-negative because $x\longmapsto\phi_k(x; t)$ is non-decreasing, while
$R_k(x; t)/d_i^{(k)}(t)$ has positive sign.
It also shows that $b_i^{(k)}(t)>0$ exactly when
the last integrand is nonzero on a Borel set $E$ with $\nu_k(E; t)>0$.
\end{proof}

It remains to prove that $\mathbf{A}$ is a nonsingular $\mathcal{M}$-matrix;
for one measure it is the identity.

\begin{theorem}\label{LFinal}
Assume the hypotheses and notation of Lemma~\ref{Ax=b}, and fix $t$. Then
\begin{align*}
\mathbf{A}(t)&=\mathbf{I}_n-\mathbf{B}(t),\\[7pt]
\mathbf{B}(t)&\ge0,\quad \rho(\mathbf{B}(t))<1.
\end{align*}
Consequently $\mathbf{A}(t)$ is a nonsingular $\mathcal{M}$-matrix and
$\mathbf{A}(t)^{-1}\ge0$. If $m\ge2$ and the support of the $k$th measure at
the fixed parameter $t$ contains more than $n_k$ points for every $k$, then
$\mathbf{B}(t)$ is irreducible and
$\mathbf{A}(t)^{-1}>0$. If
$\mathbf{u}=(u_1,\ldots,u_n)^{\mathrm T}\ge0$ and
$r\in\{1,\ldots,n\}$, the $r$th component of
$\mathbf{A}(t)^{-1}\mathbf{u}$ is positive if and only if there exist an integer $q\ge0$
and an index $s\in\{1,\ldots,n\}$ such that
\begin{align}\label{eq:reachability}
(\mathbf{B}(t)^q)_{rs}u_s>0.
\end{align}
\end{theorem}
\begin{proof}
Lemma~\ref{lemma1} gives $\mathbf{A}(t)=\mathbf{I}_n-\mathbf{B}(t)$ with
$\mathbf{B}(t)\ge0$. We first suppose that the support of the $k$th measure at
the fixed parameter $t$ has more than $n_k$ points for every $k$. For $m=1$ there is
nothing to prove. For $m\ge2$, all cross-block entries of $\mathbf{B}(t)$ are positive,
so $\mathbf{B}(t)$ is irreducible. Let $\mathbf{v}(t)>0$ be a Perron vector
satisfying
$\mathbf{B}(t)\mathbf{v}(t)=\rho(\mathbf{B}(t))\mathbf{v}(t)$. Denote
the component of $\mathbf{v}(t)$ with block index $(k,i)$ by $v_{k,i}(t)$,
and form
\begin{align}\label{eq:S-poly}
S(x; t)=\sum_{k=1}^m\sum_{i=1}^{n_k}
 v_{k,i}(t)\frac{P_n(x; t)}{x-x_{k,i}(t)}.
\end{align}
We have
$$
\frac{\partial}{\partial x}\left(\frac{S(x; t)}{P_n(x; t)}\right)
=-\sum_{k=1}^m\sum_{i=1}^{n_k}
\frac{v_{k,i}(t)}{(x-x_{k,i}(t))^2}<0.
$$
The quotient tends to opposite infinities at the two ends of every gap
between consecutive zeros of $P_n(x; t)$. It follows that $S(x; t)$ has
one simple zero in each such gap and none outside the extreme zeros. There are $m-1$ gaps between
successive blocks; hence some interval, say $\Delta_{k_0}$, contains no
zero arising from a gap between blocks. Let $H(x; t)$ be the monic polynomial
whose $n_{k_0}-1$ zeros are the zeros of $S(x; t)$ between consecutive zeros of
$Q_{n_{k_0}}(x; t)$. Interlacing gives
\begin{align}\label{eq:H-lagrange}
H(x; t)&=\sum_{i=1}^{n_{k_0}}c_i(t)
       \frac{Q_{n_{k_0}}(x; t)}{x-x_{k_0,i}(t)},\\[7pt]
 c_i(t)&=\frac{H(x_{k_0,i}(t); t)}
 {\dps\frac{\partial Q_{n_{k_0}}}{\partial x}(x_{k_0,i}(t); t)}>0.
\end{align}
In every gap between consecutive zeros of $Q_{n_{k_0}}(x; t)$, the
polynomials $H(x; t)$ and $S(x; t)$ have the same unique simple zero.
Their product therefore does not change sign there. At the zeros of
$Q_{n_{k_0}}(x; t)$,
\begin{align*}
H(x_{k_0,i}(t); t)S(x_{k_0,i}(t); t)
 &=c_i(t)v_{k_0,i}(t)
 \left(\frac{\partial Q_{n_{k_0}}}{\partial x}(x_{k_0,i}(t); t)\right)^2
 \\[7pt]
 &\hspace{12mm}{}\times
 R_{k_0}(x_{k_0,i}(t); t),
\end{align*}
Thus $H(x; t)S(x; t)$ has the constant sign of $R_{k_0}(x; t)$ on
$\Delta_{k_0}$. Its integral is nonzero: the product has only
$n_{k_0}-1$ distinct zeros there, whereas the support of the $k_0$th
measure has more than $n_{k_0}$ points.

For later use, direct substitution in \eqref{eq:S-poly} and orthogonality
give the identity
\begin{align}\label{eq:key-Av}
\int\frac{Q_{n_k}(x; t)}{x-x_{k,i}(t)}S(x; t)\,\mathrm{d}\nu_k(x; t)
 =d_i^{(k)}(t)(\mathbf{A}(t)\mathbf{v}(t))_{k,i}.
\end{align}
Indeed, the term with the same block and $j=i$ is
$v_{k,i}(t)d_i^{(k)}(t)$. If $j\ne i$ in that block, its multiplier of
$P_n(x; t)$ is
$Q_{n_k}(x; t)/((x-x_{k,i}(t))(x-x_{k,j}(t)))$, a polynomial of degree
$n_k-2$, and its integral vanishes. A term from a block $l\ne k$ is
$v_{l,j}(t)d_i^{(k)}(t)a_{i,j}^{(k,l)}(t)$. This proves
\eqref{eq:key-Av}. Combining
\eqref{eq:H-lagrange}--\eqref{eq:key-Av},
$$
\int H(x; t)S(x; t)\,\mathrm{d}\nu_{k_0}(x; t)
 =(1-\rho(\mathbf{B}(t)))\sum_{i=1}^{n_{k_0}}
 c_i(t)d_i^{(k_0)}(t)v_{k_0,i}(t).
$$
Both the integral and the sum have the sign of $R_{k_0}(x; t)$; therefore
$\rho(\mathbf{B}(t))<1$.

We now remove the assumption on the size of the supports. First,
$\mathbf{A}(t)$ is always invertible. Indeed, let
$\mathbf{w}\in\mathbb{R}^n$ satisfy
$\mathbf{A}(t)\mathbf{w}=\mathbf{0}$. Write its components as
$w_{k,i}$, and define
$$
T(x; t)=\sum_{k=1}^m\sum_{i=1}^{n_k}
w_{k,i}\frac{P_n(x; t)}{x-x_{k,i}(t)}.
$$
The same computation as in \eqref{eq:key-Av} gives
$$
\int \frac{Q_{n_k}(x; t)}{x-x_{k,i}(t)}T(x; t)
\,\mathrm{d}\nu_k(x; t)
=d_i^{(k)}(t)(\mathbf{A}(t)\mathbf{w})_{k,i}
=0,\quad 1\le i\le n_k.
$$
The displayed polynomials form a basis of $\mathbb{P}_{n_k-1}$, so $T(x; t)$
satisfies all the orthogonality conditions for the multi-index. Since
$\deg_x T(x; t)\le n-1$, normality gives $T(x; t)=0$. At a zero of
$P_n(x; t)$, however,
$$
T(x_{k,i}(t); t)=w_{k,i}
\frac{\partial P_n}{\partial x}(x_{k,i}(t); t).
$$
The zeros are simple, and therefore $w_{k,i}=0$ for every $k$ and $i$.

Choose compactly supported positive measures $\tau_k$ of infinite support in
$\Delta_k$ and replace $\mathrm{d}\nu_k(x; t)$ by
$\mathrm{d}\nu_k(x; t)+\varepsilon\,\mathrm{d}\tau_k(x)$. The preceding case gives
$\rho(\mathbf{B}_\varepsilon(t))<1$. The moments converge as
$\varepsilon\downarrow0$. Invertibility of the limiting moment system then
gives convergence of the coefficients by Cramer's rule, simplicity gives
convergence of the labelled zeros, and $d_i^{(k)}(t)\ne0$ gives convergence
of the entries of $\mathbf{B}_\varepsilon(t)$. Consequently
$\rho(\mathbf{B}(t))\le1$. Equality is impossible because the Perron root of
$\mathbf{B}(t)\ge0$ is an eigenvalue and
$\mathbf{I}_n-\mathbf{B}(t)=\mathbf{A}(t)$ is invertible. Thus
$\rho(\mathbf{B}(t))<1$ in all cases. Finally,
$$
\mathbf{A}(t)^{-1}=\sum_{q=0}^{\infty}\mathbf{B}(t)^q\ge0,
$$
which also proves the strict statement and \eqref{eq:reachability}.
\end{proof}

\begin{remark}\label{rem:saturated}
The qualification in the strict statement is necessary. If
$m=2$, $(n_1,n_2)=(1,1)$, $\nu_1=\delta_0$ and $\nu_2=\delta_2$, then
$P_2(x)=x(x-2)$, $\mathbf{A}=\mathbf{I}_2$, and both zeros are fixed under every positive
reweighting of the two atoms. More generally, if
$\#\supp\nu_k=n_k$, the zeros in the $k$th block are exactly those atoms and
the corresponding rows of $\mathbf{B}$ vanish.
\end{remark}

\begin{theorem}\label{maintheorem}
Let $I$ be an open interval, let $n_1,\ldots,n_m$ be positive integers,
and put $n=n_1+\cdots+n_m$. Let $\Delta_k=(a_k,b_k)$,
$k=1,\ldots,m$, be pairwise disjoint open intervals, and let
$\mu_1,\ldots,\mu_m$ be fixed positive Borel measures concentrated on the
respective intervals. Assume $\#\supp\mu_k\ge n_k$ for every
$k$, and let
$$
\mathrm{d}\nu_k(x; t)=\omega_k(x; t)\,\mathrm{d}\mu_k(x),\quad t\in I,
$$
where $\omega_k(x; t)$ is Borel and positive for $\mu_k$-almost every $x$.
For every $k=1,\ldots,m$ and $0\le r\le n+n_k-1$, assume that
$$
c_r^{(k)}(t)=\int x^r\omega_k(x; t)\,\mathrm{d}\mu_k(x)
$$
is finite and belongs to $C^1(I)$. Suppose that, for every $k=1,\ldots,m$
and $t\in I$, there is a Borel function
$x\longmapsto\phi_k(x; t)$ such that

\begin{align}\label{eq:main-moment-derivative}
\int(1+|x|^{n+n_k-1})|\phi_k(x; t)|
\,\mathrm{d}\nu_k(x; t)&<\infty,\\[7pt]
\frac{\mathrm{d}c_r^{(k)}}{\mathrm{d}t}(t)
 &=\int x^r\phi_k(x; t)\,\mathrm{d}\nu_k(x; t),\quad
 0\le r\le n+n_k-1.
\end{align}
When differentiation under the integral sign is justified and the displayed
absolute-integrability condition holds for the following quotient, one may
take
$$
\phi_k(x; t)=\frac{1}{\omega_k(x; t)}
\frac{\partial\omega_k}{\partial t}(x; t).
$$
For every $t\in I$, let $P_n(x; t)$ be the monic type II polynomial for
the measures $\nu_1,\ldots,\nu_m$ at parameter $t$ and the multi-index
$(n_1,\ldots,n_m)$; its existence, uniqueness and simple real zeros follow
from the support hypotheses. If, for every $k$ and $t$, the function
$x\longmapsto\phi_k(x; t)$ agrees $\nu_k$-almost everywhere at parameter
$t$ with a non-decreasing function on $\Delta_k$, then every zero of
$P_n(x; t)$ is a non-decreasing function of $t$ on $I$. If all these
representatives are non-increasing, every zero is non-increasing.
\end{theorem}
\begin{proof}
Lemma~\ref{aux1} gives $\mathbf{b}(t)\ge0$. By
Theorem~\ref{LFinal},
$\frac{\mathrm{d}\mathbf{x}}{\mathrm{d}t}(t)
=\mathbf{A}(t)^{-1}\mathbf{b}(t)\ge0$.
The decreasing case follows after changing signs. Applying the local
conclusion at every point of $I$ proves the global statement.
\end{proof}

Under the non-decreasing hypotheses of the theorem, so that
$\mathbf{b}(t)\ge0$, the proof retains the following strictness information.
It is governed by \eqref{eq:reachability} with $\mathbf{u}=\mathbf{b}(t)$. In
particular, if $m\ge2$,
$\#\supp\mu_k>n_k$ for every $k$, and one component of $\mathbf{b}(t)$ is
positive, then the derivatives of all the zeros are positive at $t$. The
decreasing case is obtained by applying the same criterion to
$-\mathbf{b}(t)$.

The proof uses only finitely many differentiated moment identities and the
$\mathcal{M}$-matrix argument.

\begin{prop}\label{prop:moving-endpoints}
Let $I$ be an open interval, let $n_1,\ldots,n_m$ be positive integers,
and put $n=n_1+\cdots+n_m$. For $k=1,\ldots,m$, let
$$
\mathrm{d}\nu_k(x; t)=\omega_k(x; t)
\mathbf 1_{(a_k(t),b_k(t))}(x)\,\mathrm{d}x,
$$
where the intervals remain pairwise disjoint, $a_k,b_k\in C^1(I)$, and
$\omega_k(x; t)>0$ in the interior. Let $P_n(x; t)$ denote the corresponding
monic type II polynomial and assume normality for every $t\in I$. For
$0\le r\le n+n_k-1$, assume that the moments are $C^1$, that
$x^r\frac{\partial\omega_k}{\partial t}(x; t)$ is integrable, and that
\begin{align*}
\frac{\mathrm{d}}{\mathrm{d}t}
\int_{a_k(t)}^{b_k(t)}x^r\omega_k(x; t)\,\mathrm{d}x
&=\int_{a_k(t)}^{b_k(t)}x^r
\frac{\partial\omega_k}{\partial t}(x; t)\,\mathrm{d}x\\[7pt]
&\quad+\frac{\mathrm{d}b_k}{\mathrm{d}t}(t)b_k(t)^r
\omega_k(b_k(t)-; t)\\[7pt]
&\quad-\frac{\mathrm{d}a_k}{\mathrm{d}t}(t)a_k(t)^r
\omega_k(a_k(t)+; t),
\end{align*}
with finite one-sided traces. If the functions
$$
 x\longmapsto\frac{1}{\omega_k(x; t)}
 \frac{\partial\omega_k}{\partial t}(x; t)
$$
are non-decreasing on $(a_k(t),b_k(t))$ and
$\frac{\mathrm{d}a_k}{\mathrm{d}t}(t),
\frac{\mathrm{d}b_k}{\mathrm{d}t}(t)\ge0$, then all zeros are
non-decreasing. If the logarithmic derivatives are non-increasing and both
endpoint derivatives are non-positive, then all zeros are non-increasing.
\end{prop}
\begin{proof}
The proof of Lemma~\ref{Ax=b} is unchanged, except that the resulting system is
$$
\mathbf{A}(t)\frac{\mathrm{d}\mathbf{x}}{\mathrm{d}t}(t)
=\widetilde{\mathbf{b}}(t),
$$
where the component of $\widetilde{\mathbf{b}}(t)$ with block index $(k,i)$
is defined by
\begin{align*}
 \widetilde b_i^{(k)}(t)d_i^{(k)}(t)
 &=\int_{a_k(t)}^{b_k(t)}
 \frac{Q_{n_k}(x; t)P_n(x; t)}{x-x_{k,i}(t)}
 \frac{\partial\omega_k}{\partial t}(x; t)\,\mathrm{d}x\\[7pt]
 &\quad+\frac{\mathrm{d}b_k}{\mathrm{d}t}(t)
 \frac{Q_{n_k}(b_k(t); t)P_n(b_k(t); t)}{b_k(t)-x_{k,i}(t)}
 \omega_k(b_k(t)-; t)\\[7pt]
 &\quad-\frac{\mathrm{d}a_k}{\mathrm{d}t}(t)
 \frac{Q_{n_k}(a_k(t); t)P_n(a_k(t); t)}{a_k(t)-x_{k,i}(t)}
 \omega_k(a_k(t)+; t).
\end{align*}
For either $e_k(t)=a_k(t)$ or $e_k(t)=b_k(t)$, cancellation gives
\begin{align*}
&\frac{Q_{n_k}(e_k(t); t)P_n(e_k(t); t)}
{d_i^{(k)}(t)(e_k(t)-x_{k,i}(t))}\\[7pt]
&\quad=\frac{R_k(e_k(t); t)}{d_i^{(k)}(t)}
(e_k(t)-x_{k,i}(t))
\left(\frac{Q_{n_k}(e_k(t); t)}
{e_k(t)-x_{k,i}(t)}\right)^2.
\end{align*}
The first factor on the right is positive. Hence the upper boundary
coefficient is positive and the lower one is negative; the latter occurs
with a minus sign in the definition of $\widetilde b_i^{(k)}(t)$. The two boundary
contributions are therefore non-negative under the stated hypotheses. The
interior term has the sign established in Lemma~\ref{aux1};
at each fixed parameter value, the static proof of
Theorem~\ref{LFinal} applies verbatim to the resulting Angelesco system and
completes the proof.
\end{proof}

In Example~\ref{Example1}, Proposition~\ref{prop:moving-endpoints} shows that
the zeros decrease with the geometric parameter $a$ when $a_1,a_2\ge0$, and
increase with $b$ when $c_1,c_2\ge0$, whenever the displayed traces and
differentiability assumptions hold.

\section{Extensions and limitations}

\subsection{A general variation identity}\label{sec:dual}

We now drop the Angelesco assumption. Fix a normal type II problem for the
multi-index $(n_1,\ldots,n_m)$, put $n=n_1+\cdots+n_m$, and suppose that
$P_n(x; t)$ has simple real zeros
$x_1(t)<\cdots<x_n(t)$. For every $i$, normality gives unique polynomials
\begin{align}\label{eq:dual-reproduction}
\deg_x\Lambda_{i,k}(x; t)&\le n_k-1,\\[7pt]
 q(x_i(t))&=\sum_{k=1}^m\int q(x)\Lambda_{i,k}(x; t)
 \,\mathrm{d}\nu_k(x; t),\quad
 q\in\mathbb{P}_{n-1}.
\end{align}
Indeed, after imposing the reproduction identity on the basis
$1,x,\ldots,x^{n-1}$, the coefficient matrix is the transpose of the
moment matrix of the normal problem. This proves existence and uniqueness of
the polynomials $\Lambda_{i,k}(x; t)$.
Define finite signed measures
\begin{align}\label{eq:sigmaik}
\mathrm{d}\sigma_{i,k}(x; t)
 =\frac{P_n(x; t)\Lambda_{i,k}(x; t)}
 {\dps\frac{\partial P_n}{\partial x}(x_i(t); t)}\,\mathrm{d}\nu_k(x; t).
\end{align}
Multiple orthogonality implies $\sigma_{i,k}(\mathbb{R}; t)=0$.

\begin{theorem}\label{thm:dual}
Under the normality and simplicity assumptions above, suppose that the
moment differentiability hypotheses of Proposition~\ref{diff} hold with
Borel admissible representing functions $x\longmapsto\phi_k(x; t)$, as
defined by \eqref{eq:representing-moments}. Then
\begin{align}\label{eq:dual-variation}
\frac{\mathrm{d}x_i}{\mathrm{d}t}(t)
 =\sum_{k=1}^m\int\phi_k(x; t)\,\mathrm{d}\sigma_{i,k}(x; t).
\end{align}
At a fixed parameter value, the inequality
$$
\sum_{k=1}^m\int f_k(x)\,\mathrm{d}\sigma_{i,k}(x; t)\ge0
$$
holds for every independent choice of bounded Borel non-decreasing functions
$x\longmapsto f_k(x)$ if and only if
\begin{align}\label{eq:separate-tails}
\sigma_{i,k}((y,\infty); t)\ge0,\quad
k=1,\ldots,m,\quad y\in\mathbb{R}.
\end{align}
On the diagonal, the inequality
$$
\sum_{k=1}^m\int f(x)\,\mathrm{d}\sigma_{i,k}(x; t)\ge0
$$
holds for every bounded Borel non-decreasing function $x\longmapsto f(x)$
if and only if
\begin{align}\label{eq:combined-tails}
\sum_{k=1}^m\sigma_{i,k}((y,\infty); t)\ge0,\quad y\in\mathbb{R}.
\end{align}
The same equivalences hold for unbounded Borel functions satisfying
$$
\sum_{k=1}^m\int|f_k(x)|\,\mathrm{d}|\sigma_{i,k}|(x; t)<\infty
$$
in the independent case, and the analogous condition with $f_1=\cdots=f_m=f$
in the diagonal case. Under \eqref{eq:separate-tails}, independent
non-increasing admissible representing functions give non-positive
derivatives of the zeros. Under \eqref{eq:combined-tails}, the same
conclusion holds when the representing functions are common, that is,
$\phi_1(x; t)=\cdots=\phi_m(x; t)$ for every $x$.
\end{theorem}
\begin{proof}
Differentiate first the monomial orthogonality conditions. At the fixed
parameter value, take their linear combination with the coefficients of
$\Lambda_{i,k}(x; t)$; these coefficients are held fixed in this operation.
The products that occur have degree at most $n+n_k-1$, so the representing
identities in Proposition~\ref{diff} apply. Differentiation with respect to
$t$ leaves the leading coefficient unchanged. Hence, as a polynomial in
$x$, the function $\frac{\partial P_n}{\partial t}(x; t)$ belongs to
$\mathbb{P}_{n-1}$, and \eqref{eq:dual-reproduction} gives
$$
\frac{\partial P_n}{\partial t}(x_i(t); t)
 =-\sum_{k=1}^m\int P_n(x; t)\Lambda_{i,k}(x; t)\phi_k(x; t)
 \,\mathrm{d}\nu_k(x; t).
$$
Together with
$$
 \frac{\mathrm{d}x_i}{\mathrm{d}t}(t)
 =-\frac{\dps\frac{\partial P_n}{\partial t}(x_i(t); t)}
 {\dps\frac{\partial P_n}{\partial x}(x_i(t); t)},
$$
this proves
\eqref{eq:dual-variation}.

It remains to justify the equivalences, including functions with jumps. Let
$\sigma$ be a finite signed measure of mass zero and suppose its open right
tails are non-negative. After subtracting a lower bound, every bounded
non-decreasing function $x\longmapsto\phi(x)$ may be assumed non-negative
and has the representation
$$
\phi(x)=\int_0^\infty\mathbf 1_{\{\phi(x)>s\}}\,\mathrm{d}s.
$$
Each superlevel set is empty, all of $\mathbb{R}$, $(c,\infty)$, or
$[c,\infty)$. Closed tails are non-negative as well, since
$\sigma([c,\infty))=\lim_{y\uparrow c}\sigma((y,\infty))$. Fubini's theorem
therefore yields $\int\phi(x)\,\mathrm{d}\sigma(x)\ge0$. Truncation and dominated
convergence give the integrable case. Conversely, testing with
$\phi(x)=\mathbf 1_{(y,\infty)}(x)$ gives each tail inequality. Apply this argument
separately to the $\sigma_{i,k}$, or to their sum, respectively.
\end{proof}

A negative separate or combined tail is realised by multiplying,
respectively, the corresponding measure or all the measures by
$\vartheta(x; s)=\exp\{s\mathbf 1_{(y,\infty)}(x)\}$. At $s=0$ its
logarithmic derivative is the bounded non-decreasing function
$\mathbf 1_{(y,\infty)}(x)$, so
\eqref{eq:dual-variation} makes the derivative of the relevant zero negative.
Normality and simplicity persist for sufficiently small $|s|$ by continuity,
as in Proposition~\ref{diff}. For Angelesco systems,
Theorem~\ref{maintheorem} rules out such a negative separate tail; hence
\eqref{eq:separate-tails} always holds. For a common deformation, the
combined tails are decisive.

The extension of Theorem~\ref{maintheorem} obtained by imposing Markov's
condition on each measure of an AT system is false, even in degree two.

\begin{eje}\label{ex:AT}
On $(0,1)$ take the multi-index $(1,1)$ and
\begin{align}\label{eq:AT-weights}
\mathrm{d}\nu_1(x; t)=\omega_1(x; t)\,\mathrm{d}x
 &=x^{6/5}e^{tx^2}\,\mathrm{d}x,\\[7pt]
 \mathrm{d}\nu_2(x; t)=\omega_2(x; t)\,\mathrm{d}x
 &=x^{5/4}\,\mathrm{d}x.
\end{align}
For $|t|<1/40$ this is an AT system at $(1,1)$, since
$$
\frac{\omega_1(x; t)}{\omega_2(x; t)}
 \frac{\partial}{\partial x}
 \left(\frac{\omega_2(x; t)}{\omega_1(x; t)}\right)
 =\frac{1}{20x}-2tx>0.
$$
At $t=0$ the two logarithmic derivatives are
\begin{align*}
\frac{1}{\omega_1(x; 0)}
\frac{\partial\omega_1}{\partial t}(x; 0)=x^2,\quad
\frac{1}{\omega_2(x; 0)}
\frac{\partial\omega_2}{\partial t}(x; 0)=0.
\end{align*}
Both functions are non-decreasing. Nevertheless, the two zeros satisfy

\begin{align}\label{eq:AT-negative}
\frac{\mathrm{d}x_1}{\mathrm{d}t}(0)&=-\frac{8800}{8277}<0,\\[7pt]
\frac{\mathrm{d}x_2}{\mathrm{d}t}(0)&=-\frac{950400}{2298247}<0.
\end{align}
\end{eje}
\begin{proof}
Write $P_{(1,1)}(x; t)=x^2+a(t)x+b(t)$. Since
$\int_0^1x^\alpha\,\mathrm{d}x=(\alpha+1)^{-1}$, the two orthogonality
conditions at $t=0$ reduce to

\begin{align*}
\frac{1}{21}+\frac{a(0)}{16}+\frac{b(0)}{11}=0,\quad
\frac{1}{17}+\frac{a(0)}{13}+\frac{b(0)}{9}=0.
\end{align*}
The determinant of the system for $a(0)$ and $b(0)$ is
$-1/20592$, and therefore
\begin{align}\label{eq:AT-system}
P_{(1,1)}(x; 0)
 =x^2-\frac{416}{357}x+\frac{33}{119}
 =\left(x-\frac{1}{3}\right)\left(x-\frac{99}{119}\right).
\end{align}
Because
$$
 \frac{\partial P_{(1,1)}}{\partial t}(x; 0)
 =\frac{\mathrm{d}a}{\mathrm{d}t}(0)x+\frac{\mathrm{d}b}{\mathrm{d}t}(0),
$$
differentiation of the
orthogonality conditions gives
\begin{align*}
\int_0^1x^{6/5}\left(\frac{\partial P_{(1,1)}}{\partial t}(x; 0)
 +x^2P_{(1,1)}(x; 0)\right)\,\mathrm{d}x&=0,\\[7pt]
 \int_0^1x^{5/4}\frac{\partial P_{(1,1)}}{\partial t}(x; 0)\,\mathrm{d}x&=0.
\end{align*}
On evaluating the integrals, these equations become
\begin{align*}
\frac{1}{16}\frac{\mathrm{d}a}{\mathrm{d}t}(0)
+\frac{1}{11}\frac{\mathrm{d}b}{\mathrm{d}t}(0)
=-\frac{50}{77469},\quad
\frac{1}{13}\frac{\mathrm{d}a}{\mathrm{d}t}(0)
+\frac{1}{9}\frac{\mathrm{d}b}{\mathrm{d}t}(0)=0.
\end{align*}
The same nonzero determinant gives
\begin{align*}
\frac{\mathrm{d}a}{\mathrm{d}t}(0)=\frac{114400}{77469},\quad
\frac{\mathrm{d}b}{\mathrm{d}t}(0)=-\frac{26400}{25823}.
\end{align*}
Hence
\begin{align}\label{eq:AT-root-derivative}
\frac{\partial P_{(1,1)}}{\partial t}(x; 0)
 =\frac{114400}{77469}\left(x-\frac{9}{13}\right).
\end{align}
Now $1/3<9/13<99/119$. Thus $\frac{\partial P_{(1,1)}}{\partial t}(x; 0)$ is
negative at the first zero and positive at the second, whereas
$\frac{\partial P_{(1,1)}}{\partial x}(x; 0)$ has the same signs there.
Differentiating $P_{(1,1)}(x_i(t); t)=0$ gives
$$
 \frac{\mathrm{d}x_i}{\mathrm{d}t}(0)
 =-\frac{\dps\frac{\partial P_{(1,1)}}{\partial t}(x_i(0); 0)}
 {\dps\frac{\partial P_{(1,1)}}{\partial x}(x_i(0); 0)},
$$
and direct substitution gives
\eqref{eq:AT-negative}.
\end{proof}

Thus the componentwise Markov condition fails on a common support, and
Theorem~\ref{thm:dual} gives its replacement. The following theorem records a
case in which the combined condition is automatic.

\begin{theorem}\label{thm:minimal-support}
Let $u_1,\ldots,u_n$ be a Chebyshev system on an interval containing
$y_0<\cdots<y_n$, and let
\begin{align}\label{eq:minimal-measure}
\mathrm{d}\nu(x; t)=\sum_{j=0}^n a_j
 e^{t\phi(y_j)}\,\mathrm{d}\delta_{y_j}(x),\quad
 a_j>0,\quad t\in\mathbb{R},
\end{align}
where the values $\phi(y_j)$ are non-decreasing with $j$. If $P_n(x; t)$ is
the monic polynomial determined by
$$
\int P_n(x; t)u_k(x)\,\mathrm{d}\nu(x; t)=0,\quad k=1,\ldots,n,
$$
then all its zeros are non-decreasing functions of $t$. They are strictly
increasing if the values $\phi(y_0),\ldots,\phi(y_n)$ are not all equal.
\end{theorem}
\begin{proof}
The evaluation determinants of a Chebyshev system on ordered points have one
common nonzero sign. After multiplying one of the functions by $-1$, if
necessary, assume that this sign is positive, and put
$D(x)=\prod_{j=0}^n(x-y_j)$. Expanding the null vector of
the $n$ by $n+1$ evaluation matrix by cofactors gives
\begin{align}\label{eq:minimal-fraction}
\frac{P_n(x; t)}{D(x)}
 &=\sum_{j=0}^n\frac{r_j(t)}{x-y_j},\\[7pt]
 r_j(t)&=\frac{c_j e^{-t\phi(y_j)}}
 {\dps\sum_{l=0}^n c_l e^{-t\phi(y_l)}},
\end{align}
where
$$
c_j=\frac{\det[u_k(y_l)]_{\substack{1\le k\le n\\[7pt]0\le l\le n,\ l\ne j}}}
 {a_j\left|\dps\frac{\mathrm{d}D}{\mathrm{d}x}(y_j)\right|}>0.
$$
Let $\Delta_j$ denote the determinant in the numerator. Cofactor expansion
gives, for one constant $\lambda(t)$,
$$
a_j e^{t\phi(y_j)}P_n(y_j; t)=\lambda(t)(-1)^j\Delta_j.
$$
At $y_j$, the derivative of $D$ has sign $(-1)^{n-j}$, so all residues
have the same sign. Moreover,
$$
\frac{P_n(x; t)}{D(x)}=x^{-1}+O(x^{-2})
$$
at infinity, and therefore the residues sum to one. This fixes $\lambda(t)$
and proves \eqref{eq:minimal-fraction} with $r_j(t)>0$.

Write the right-hand side of \eqref{eq:minimal-fraction} as $F(x; t)$, and
observe first that every $r_j(t)$ is positive. On each interval
$(y_j,y_{j+1})$, the function $x\longmapsto F(x; t)$ decreases from $+\infty$ to
$-\infty$, because
$$
 \frac{\partial F}{\partial x}(x; t)
 =-\sum_{j=0}^n\frac{r_j(t)}{(x-y_j)^2}<0.
$$
It follows that $P_n(x; t)$ has exactly one simple zero in each of these
intervals. Let $\xi(t)$ be any one of them. Extend the values $\phi(y_j)$
to a non-decreasing function on the interval, and put
$$
\overline\phi(t)=\sum_{l=0}^n r_l(t)\phi(y_l).
$$
Direct differentiation of the second identity in
\eqref{eq:minimal-fraction} gives
$$
\frac{\mathrm{d}r_j}{\mathrm{d}t}(t)
=r_j(t)\bigl(\overline\phi(t)-\phi(y_j)\bigr).
$$
At a zero $\xi(t)$, the preceding formula gives
$\frac{\partial F}{\partial x}(\xi(t); t)<0$, while
$F(\xi(t); t)=0$ gives
\begin{align*}
\frac{\partial F}{\partial t}(\xi(t); t)
&=-\sum_{j=0}^n r_j(t)
\frac{\phi(y_j)-\overline\phi(t)}{\xi(t)-y_j}\\[7pt]
&=-\sum_{j=0}^n r_j(t)
\frac{\phi(y_j)-\phi(\xi(t))}{\xi(t)-y_j}\ge0.
\end{align*}
The middle expression may be changed into the last one because their
difference is a constant multiple of $F(\xi(t); t)$. Every quotient in the
last sum is non-positive. If the values $\phi(y_j)$
are not all equal, at least one of these quotients is strictly negative for
every zero $\xi(t)$. Hence
$$
\frac{\mathrm{d}\xi}{\mathrm{d}t}(t)
 =-\frac{\dps\frac{\partial F}{\partial t}(\xi(t); t)}
 {\dps\frac{\partial F}{\partial x}(\xi(t); t)}\ge0.
$$
The last inequality is strict when the values $\phi(y_j)$ are not all equal.
\end{proof}

\subsection{Nikishin systems}\label{sec:nikishin}

We use the convention
$$
\widehat\tau(z)=\int\frac{1}{z-y}\,\mathrm{d}\tau(y).
$$
Let $\sigma_1,\ldots,\sigma_m$ be finite positive Borel measures with
infinite compact supports. Write $\Delta_j$ for the convex hull of
$\supp\sigma_j$, and assume that consecutive intervals $\Delta_j$ are
disjoint. Define
\begin{align*}
s_{j,j}&=\sigma_j,\\[7pt]
\mathrm{d}s_{j,k}(x)&=\widehat s_{j+1,k}(x)\,\mathrm{d}\sigma_j(x),
\quad j<k.
\end{align*}
The ordered system of measures
$$
\mathcal{N}(\sigma_1,\ldots,\sigma_m)
=(s_{1,1},\ldots,s_{1,m})
$$
is a Nikishin system. Fix a multi-index $(n_1,\ldots,n_m)$ and put
$n=n_1+\cdots+n_m$. Write its monic type II polynomial as
$P_n$, and choose constants
$\varepsilon_1=1$, $\varepsilon_k\in\{-1,1\}$ so that
$\nu_k=\varepsilon_k s_{1,k}$ is positive. This does not alter orthogonality.
Such systems are perfect \cite[Theorem~1.2]{FidalgoPerfect}; their type II
zeros are simple and interior to $\Delta_1$, and neighbouring polynomials
strictly interlace \cite[Corollaries~4.5 and~5.5]{KozhanVaktnasZeros}.

We first deform the first generator. Keep $\sigma_2,\ldots,\sigma_m$ fixed and set
$\mathrm{d}\sigma_1(x; t)=\omega(x; t)\,\mathrm{d}\mu(x)$. Define
$\mathrm{d}s_{1,1}^0(x)=\mathrm{d}\mu(x)$ and
$\mathrm{d}s_{1,k}^0(x)=\widehat s_{2,k}(x)\,\mathrm{d}\mu(x)$ for $k\ge2$.
Then, for every $k$,
\begin{align}\label{eq:nik-common}
\mathrm{d}s_{1,k}(x; t)=\omega(x; t)\,\mathrm{d}s_{1,k}^0(x).
\end{align}
Thus every positive multiplicative deformation of the first generator has
the same logarithmic derivative on all components. With the dual polynomials
from \eqref{eq:dual-reproduction}, write, for $x\in\Delta_1$,
\begin{align*}
w_1(x)&=1,\quad
w_k(x)=\varepsilon_k\,\widehat s_{2,k}(x)>0,\quad
k=2,\ldots,m,\\[7pt]
\mathcal{K}_i(x; t)&=\sum_{k=1}^m\Lambda_{i,k}(x; t)w_k(x).
\end{align*}
The combined measure in \eqref{eq:combined-tails} is therefore
\begin{align}\label{eq:nik-dual}
\mathrm{d}\Sigma_i(x; t)&=
 \frac{P_n(x; t)\mathcal{K}_i(x; t)}
 {\dps\frac{\partial P_n}{\partial x}(x_i(t); t)}
 \,\mathrm{d}\sigma_1(x; t),\\[7pt]
 \Sigma_i(\mathbb{R}; t)&=0.
\end{align}
Assume that the induced positive component measures satisfy the moment
differentiability hypotheses of Proposition~\ref{diff}, with the common
admissible representing function
$$
\phi(x; t)=\frac{1}{\omega(x; t)}
 \frac{\partial\omega}{\partial t}(x; t).
$$
Then
\begin{align}\label{eq:nik-variation}
\frac{\mathrm{d}x_i}{\mathrm{d}t}(t)
 =\int\phi(x; t)\,\mathrm{d}\Sigma_i(x; t).
\end{align}
By Theorem~\ref{thm:dual}, the Markov conclusion for every common
non-decreasing logarithmic derivative is equivalent to
\begin{align}\label{eq:nik-tail}
\Sigma_i((y,\infty); t)\ge0,\quad y\in\mathbb{R}.
\end{align}
Perfectness and interlacing alone do not prove these tails: the multiple
Christoffel--Darboux kernel is generally non-symmetric and has no scalar
sum-of-squares representation \cite{DaemsKuijlaars,SwiderskiVanAssche}.

We next consider an exterior Christoffel factor. Write $\Delta_1=[a,b]$. Given
$z\in\mathbb{R}\setminus[a,b]$, let
$\widetilde P_n^{[z]}$ be the type II polynomial after replacing
$\sigma_1$ by $|x-z|\sigma_1$.

\begin{theorem}\label{thm:nik-christoffel}
For every multi-index $(n_1,\ldots,n_m)$, with $n=n_1+\cdots+n_m$, and every
$z\notin[a,b]$, if
$x_1<\cdots<x_n$ and
$\widetilde x_1^{[z]}<\cdots<\widetilde x_n^{[z]}$ are the zeros of $P_n$
and $\widetilde P_n^{[z]}$, respectively, then
\begin{align*}
 z<a&:\quad x_1<\widetilde x_1^{[z]}<x_2<\cdots
       <x_n<\widetilde x_n^{[z]},\\[7pt]
 z>b&:\quad \widetilde x_1^{[z]}<x_1<\cdots
       <\widetilde x_n^{[z]}<x_n.
\end{align*}
The same conclusion holds for an AT system on $[a,b]$ multiplied by the common factor
$|x-z|$, provided the multi-indices $(n_1,\ldots,n_m)$ and
$(n_1,\ldots,n_j+1,\ldots,n_m)$ are normal for some $j$, their type II
polynomials strictly interlace, and the transformed problem at
$(n_1,\ldots,n_m)$ is normal.
\end{theorem}
\begin{proof}
Choose $j\in\{1,\ldots,m\}$; in the AT case choose one satisfying the stated
hypotheses. Put $p(x)=P_n(x)$, and let $q$ be the monic type II polynomial
associated with the multi-index $(n_1,\ldots,n_j+1,\ldots,n_m)$. The
common-factor relation and the orthogonality relations give
\begin{align}\label{eq:christoffel-formula}
\widetilde P_n^{[z]}(x)
 =\frac{q(x)-q(z)p(x)/p(z)}{x-z}.
\end{align}
Indeed, $|x-z|$ is a constant nonzero multiple of $x-z$ on $[a,b]$.
The numerator vanishes at $z$. For every component measure and every
required exponent $r$, its transformed orthogonality follows from
\begin{align*}
&\int x^r\widetilde P_n^{[z]}(x)|x-z|
\,\mathrm{d}\nu_l(x)\\[7pt]
&\quad=\pm\int x^r q(x)\,\mathrm{d}\nu_l(x)
\mp\frac{q(z)}{p(z)}\int x^r p(x)\,\mathrm{d}\nu_l(x)=0.
\end{align*}
The quotient in \eqref{eq:christoffel-formula} is monic of degree $n$;
normality therefore identifies it with the transformed polynomial. Set
$R(x)=q(x)/p(x)$. Strict
interlacing and monicity imply
\begin{align}\label{eq:R-fraction}
 \frac{q(x)}{p(x)}
 &=x-\gamma+\sum_{i=1}^n\frac{r_i}{x-x_i},\\[7pt]
 r_i&=\frac{q(x_i)}{\dps\frac{\mathrm{d}p}{\mathrm{d}x}(x_i)}<0.
\end{align}
The last sign follows because $q(x_i)$ and
$\frac{\mathrm{d}p}{\mathrm{d}x}(x_i)$ have opposite signs at every
interlacing zero.
Consequently
$$
 \frac{\mathrm{d}R}{\mathrm{d}x}(x)
 =1-\sum_{i=1}^n\frac{r_i}{(x-x_i)^2}>0
$$
on every component of $\mathbb{R}\setminus\{x_1,\ldots,x_n\}$. On each such
component $R(x)$ increases from $-\infty$ to $+\infty$. Hence the equation
$R(x)=R(z)$ has one solution in every component; the solution in the component
containing $z$ is $z$ itself. The remaining solutions are precisely the
transformed zeros and give the displayed order. The Nikishin hypotheses were
recalled above; related Christoffel interlacings appear in
\cite[Corollary~3.23]{KozhanVaktnasChristoffel}.
\end{proof}

Multiplication by the positive common exterior factor preserves the AT
property. For an AT system on a closed interval, normality and simple interior type II
zeros are automatic at every multi-index for which the system is AT; if it is
AT at the two indicated neighbouring multi-indices, strict interlacing is
automatic as well
\cite[Theorem~3.6 and Corollaries~4.4 and~5.4]{KozhanVaktnasZeros}.
Jacobi--Pi\~neiro systems are classical examples on a compact interval, while
multiple Laguerre systems give standard unbounded examples
\cite[Sections~23.3.2 and~23.4]{I05}.

Iteration treats products of exterior factors: points on the left move all
zeros to the right, points on the right move them to the left, and mixed
products are monotone in each exterior coordinate. No density is needed. For
the continuous form, let
$R(x)=q(x)/p(x)$ refer to the fixed base system generated by $\mu$, and let
$z:I\longrightarrow\mathbb{R}\setminus[a,b]$ be a $C^1$ function whose values
remain in one component. If
$$
\mathrm{d}\sigma_1(x; t)=|x-z(t)|\,\mathrm{d}\mu(x),
$$
denote the transformed zeros by $\xi_i(z(t))$. As functions of the exterior
point, they satisfy
\begin{align}\label{eq:moving-pole}
\frac{\mathrm{d}\xi_i}{\mathrm{d}z}(z)
 =\frac{\dps\frac{\mathrm{d}R}{\mathrm{d}x}(z)}
 {\dps\frac{\mathrm{d}R}{\mathrm{d}x}(\xi_i(z))}>0.
\end{align}
This follows by differentiating $R(\xi_i(z))=R(z)$. Moreover,
$$
 \frac{1}{|x-z(t)|}\frac{\partial}{\partial t}|x-z(t)|
 =-\frac{1}{x-z(t)}\frac{\mathrm{d}z}{\mathrm{d}t}(t).
$$
This function is non-decreasing with respect to $x$ precisely when
$\frac{\mathrm{d}z}{\mathrm{d}t}(t)\ge0$.

The preceding result yields a larger class of deformations. Fix
$x_0\in(a,b)$ and put
$$
\kappa_z(x)=-\frac{1}{x-z}+\frac{1}{x_0-z}.
$$
\begin{theorem}\label{thm:resolvent}
For the standard infinite-support Nikishin system above, fix a multi-index
$(n_1,\ldots,n_m)$, put $n=n_1+\cdots+n_m$, and let $I$ be an open interval.
Let $x_1(t),\ldots,x_n(t)$ be the zeros of $P_n(x; t)$ and,
for this multi-index, let the signed measures $\Sigma_i$ be those defined by
\eqref{eq:nik-dual}.
Assume that, for every $t\in I$, the induced positive component measures
satisfy the moment differentiability hypotheses of Proposition~\ref{diff}
with the common admissible representing function
$$
\phi(x; t)=\frac{1}{\omega(x; t)}
\frac{\partial\omega}{\partial t}(x; t).
$$
Suppose that, for every $t\in I$, this function satisfies

\begin{align}\label{eq:resolvent-cone}
\phi(x; t)=C(t)+\alpha(t)(x-x_0)
 +\int_{\mathbb{R}\setminus[a,b]}\kappa_z(x)\,\mathrm{d}\eta_t(z),
\end{align}
where $C(t)\in\mathbb{R}$, $\alpha(t)\ge0$, $\eta_t$ is a positive Borel
measure, the displayed integral is well-defined, and
\begin{align}\label{eq:resolvent-fubini}
\int_{\mathbb{R}\setminus[a,b]}\int_{[a,b]}|\kappa_z(x)|
\,\mathrm{d}|\Sigma_i|(x; t)\,\mathrm{d}\eta_t(z)<\infty
\end{align}
for $i=1,\ldots,n$ and every $t\in I$. Then
$\frac{\mathrm{d}x_i}{\mathrm{d}t}(t)\ge0$ for every zero of
$P_n(x; t)$, with strict inequality when $\eta_t\ne0$. Thus the zeros are
non-decreasing, and strictly increasing on intervals where $\eta_t\ne0$
throughout. The signs reverse for the negative of a function in
\eqref{eq:resolvent-cone}.
\end{theorem}
\begin{proof}
Fix $t\in I$ and $z\notin[a,b]$, and write the current first generator as
$\mathrm{d}\sigma_1(x; t)=|x-z|\,\mathrm{d}\tau_z(x)$. If the point $z$ is replaced by $z+u$
and $\tau_z$ is kept fixed, its logarithmic derivative at $u=0$ is
$-\frac{1}{x-z}$. Formula~\eqref{eq:moving-pole} therefore shows that the right-hand
side of \eqref{eq:nik-variation} is strictly positive for this function. It
has the same value for $\kappa_z$, because the added term is constant in $x$
and $\Sigma_i(\mathbb{R}; t)=0$. Moreover,
$$
L^2\kappa_{-L}(x)\longrightarrow x-x_0
$$
uniformly on $[a,b]$. Thus the right-hand side of
\eqref{eq:nik-variation} is non-negative on the affine term. The constant
term gives zero, and \eqref{eq:resolvent-fubini} permits the use of Fubini's
theorem for the integral term. Since its integrand is strictly positive for
every exterior $z$, the result is strict whenever $\eta_t$ is nonzero.
\end{proof}

Since the multi-index was arbitrary, the theorem applies simultaneously to
every multi-index for which its moment and integrability hypotheses hold.

For finite supports, require, for every $t\in I$, the following conditions for
$\eta_t$-almost every $z$ and, when $\alpha(t)>0$, along $z=-L$ as
$L\longrightarrow\infty$: there is an index $j=j(t,z)\in\{1,\ldots,m\}$
such that the auxiliary system obtained by division by $|x-z|$ is normal at
$(n_1,\ldots,n_m)$ and $(n_1,\ldots,n_j+1,\ldots,n_m)$, the corresponding
type II polynomials strictly interlace, and the restored problem is normal at
$(n_1,\ldots,n_m)$. The same conditions give the common AT version.

Taking $\phi(x; t)=x$ shows that the zeros for
$e^{tx}\,\mathrm{d}\mu$ are non-decreasing. For $c<a$,
\begin{align}\label{eq:left-log}
\log(x-c)-\log(x_0-c)
 =\int_{-\infty}^{c}\kappa_z(x)\,\mathrm{d}z,
\end{align}
so the zeros for $(x-c)^t\,\mathrm{d}\mu$ are strictly increasing.
For $c>b$,
\begin{align}\label{eq:right-log}
-\log(c-x)+\log(c-x_0)=\int_c^\infty\kappa_z(x)\,\mathrm{d}z,
\end{align}
so the zeros for $(c-x)^t\,\mathrm{d}\mu$ are strictly decreasing. Products
on either fixed side follow by iteration.

We finally turn to an internal generator, for which the conclusion may fail.
For $m=2$, keep $\sigma_1$ fixed. If
$\mathrm{d}\sigma_2(y; t)=\eta(y; t)\,\mathrm{d}\rho(y)$ and
$$
 \psi(y; t)=\frac{1}{\eta(y; t)}
 \frac{\partial\eta}{\partial t}(y; t),
$$
the induced logarithmic derivatives on $\Delta_1$ are
\begin{align}\label{eq:internal-logarithmic-derivative}
\phi_1(x; t)&=0,\\[7pt]
 \phi_2(x; t)&=
 \frac{\dps\int \dps\frac{\psi(y; t)}{x-y}\,
 \mathrm{d}\sigma_2(y; t)}
 {\dps\int \dps\frac{1}{x-y}\,\mathrm{d}\sigma_2(y; t)}.
\end{align}
Whenever differentiation under the integral sign is justified, normalise the
constant-sign Cauchy kernel to the probability measure
$$
\mathrm{d}\pi_{x,t}(y)=
\frac{\dps\frac{1}{x-y}\,\mathrm{d}\sigma_2(y; t)}
{\dps\int \dps\frac{1}{x-v}\,\mathrm{d}\sigma_2(v; t)}.
$$
Then
\begin{align}\label{eq:internal-covariance}
\frac{\partial\phi_2}{\partial x}(x; t)
 &=\frac{1}{2}\iint
 \bigl(\psi(y; t)-\psi(v; t)\bigr)
 \left(-\frac{1}{x-y}+\frac{1}{x-v}\right)
 \,\mathrm{d}\pi_{x,t}(y)\,\mathrm{d}\pi_{x,t}(v)\le0
\end{align}
when the function $y\longmapsto\psi(y; t)$ is non-decreasing.
Indeed, $y\longmapsto-1/(x-y)$ is strictly decreasing on $\Delta_2$, so
the two factors in the integrand have opposite signs.
\begin{eje}\label{ex:internal-nikishin}
Let
\begin{align*}
\mathrm{d}\sigma_1(x)&=2\,\mathrm{d}\delta_1(x)
 +5\,\mathrm{d}\delta_4(x)+10\,\mathrm{d}\delta_7(x),\\[7pt]
\mathrm{d}\sigma_2(y; t)&=6\,\mathrm{d}\delta_{15}(y)
 +e^t\,\mathrm{d}\delta_{22}(y).
\end{align*}
Thus
$\mathrm{d}\sigma_2(y; t)=e^{t\psi(y)}\,\mathrm{d}\sigma_2(y; 0)$, where
$\psi(15)=0<1=\psi(22)$. The positive normalisation of the second measure
of this finite-support Nikishin construction is
\begin{align}\label{eq:internal-example-measure}
 \mathrm{d}\nu_2(x; t)&=-\mathrm{d}s_{1,2}(x; t)
 =w(x; t)\,\mathrm{d}\sigma_1(x),\\[7pt]
 w(x; t)&=\frac{6}{15-x}+\frac{e^t}{22-x}.
\end{align}
For the multi-index $(1,1)$, both zeros of the corresponding type II
polynomial are decreasing at $t=0$.
\end{eje}
\begin{proof}
Write $P_{(1,1)}(x; t)=x^2+a(t)x+b(t)$, and denote the three support
points of $\sigma_1$ by $\zeta_1=1$, $\zeta_2=4$, $\zeta_3=7$, with
respective masses $m_1=2$, $m_2=5$, $m_3=10$. Direct summation gives
\begin{align*}
&\sum_{l=1}^3m_l=17,\quad 
\ \ \ \ \sum_{l=1}^3m_lw(\zeta_l; 0)=\frac{8402}{693},\\[7pt]
&\sum_{l=1}^3m_l\zeta_l=92,\quad 
\ \ \sum_{l=1}^3m_l\zeta_lw(\zeta_l; 0)=\frac{97213}{1386},\\[7pt]
&\sum_{l=1}^3m_l\zeta_l^2=572,\quad
\sum_{l=1}^3m_l\zeta_l^2w(\zeta_l; 0)=\frac{622591}{1386}.
\end{align*}
Thus, at $t=0$, the two orthogonality conditions are
\begin{align*}
 17b(0)+92a(0)=-572,\quad
 \frac{8402}{693}b(0)+\frac{97213}{1386}a(0)
 =-\frac{622591}{1386}.
\end{align*}
The determinant of this system is $35551/462>0$; hence the problem is normal
and its coefficients depend differentiably on $t$ near zero. Solving gives
\begin{align}\label{eq:internal-example-polynomial}
 P_{(1,1)}(x; 0)=x^2-\frac{324053}{35551}x
 +\frac{557512}{35551}.
\end{align}
Since $\frac{\partial w}{\partial t}(x; 0)=\frac{1}{22-x}$, direct
substitution in \eqref{eq:internal-example-polynomial} also gives
$$
\sum_{l=1}^3\frac{m_lP_{(1,1)}(\zeta_l; 0)}{22-\zeta_l}
=-\frac{810}{35551}.
$$
Differentiating the two orthogonality conditions therefore gives
\begin{align*}
 17\frac{\mathrm{d}b}{\mathrm{d}t}(0)
 +92\frac{\mathrm{d}a}{\mathrm{d}t}(0)=0,\quad
 \frac{8402}{693}\frac{\mathrm{d}b}{\mathrm{d}t}(0)
 +\frac{97213}{1386}\frac{\mathrm{d}a}{\mathrm{d}t}(0)
 =\frac{810}{35551}.
\end{align*}
Consequently
\begin{align}\label{eq:internal-example-derivative}
 \frac{\partial P_{(1,1)}}{\partial t}(x; 0)
 =\frac{6361740}{35551^2}\left(x-\frac{92}{17}\right).
\end{align}
Moreover,
\begin{align*}
 P_{(1,1)}(2; 0)=\frac{51610}{35551}>0,\quad 
 P_{(1,1)}(3; 0)=-\frac{94688}{35551}<0,\\[7pt]
 P_{(1,1)}(6; 0)=-\frac{106970}{35551}<0,\quad 
 P_{(1,1)}(7; 0)=\frac{31140}{35551}>0.
\end{align*}
Thus the two simple zeros satisfy
$2<x_1(0)<3<92/17<6<x_2(0)<7$. It follows from
\eqref{eq:internal-example-derivative} that the parameter derivative of the
polynomial is negative at $x_1(0)$ and positive at $x_2(0)$, whereas its
derivative with respect to $x$ has the same signs. Therefore
$$
 \frac{\mathrm{d}x_i}{\mathrm{d}t}(0)
 =-\frac{\dps\frac{\partial P_{(1,1)}}{\partial t}(x_i(0); 0)}
 {\dps\frac{\partial P_{(1,1)}}{\partial x}(x_i(0); 0)}<0,\quad i=1,2.
$$
The same obstruction occurs for generators of infinite support. For
$\varepsilon>0$, set
\begin{align*}
\mathrm{d}\sigma_{1,\varepsilon}(x)
&=2\,\mathrm{d}\delta_1(x)+5\,\mathrm{d}\delta_4(x)
+10\,\mathrm{d}\delta_7(x)
+\varepsilon\mathbf 1_{[1,7]}(x)\,\mathrm{d}x,\\[7pt]
\mathrm{d}\sigma_{2,\varepsilon}(y; t)
&=6\,\mathrm{d}\delta_{15}(y)+e^t\,\mathrm{d}\delta_{22}(y)
+\varepsilon e^{t(y-15)/7}\mathbf 1_{[15,22]}(y)\,\mathrm{d}y.
\end{align*}
Both generators have infinite support on disjoint compact intervals and, with
$\psi(y)=(y-15)/7$ on $[15,22]$, satisfy
$\mathrm{d}\sigma_{2,\varepsilon}(y; t)
=e^{t\psi(y)}\,\mathrm{d}\sigma_{2,\varepsilon}(y; 0)$.

Let $P_{(1,1),\varepsilon}(x; t)$ and
$x_{1,\varepsilon}(t),x_{2,\varepsilon}(t)$ be the resulting polynomial and
zeros. The moment system and its derivative with respect to $t$ depend
continuously on $\varepsilon$; at $\varepsilon=0$ its determinant is
$35551/462$ and the
zeros are simple. Cramer's rule and implicit differentiation give
$$
\lim_{\varepsilon\downarrow0}
\frac{\mathrm{d}x_{i,\varepsilon}}{\mathrm{d}t}(0)
=\frac{\mathrm{d}x_i}{\mathrm{d}t}(0)<0,\quad i=1,2.
$$
The two derivatives remain negative for every sufficiently small
$\varepsilon>0$. This proves the assertion within the standard
infinite-support class.
\end{proof}

Thus monotonicity of the logarithmic derivative of an internal generator does
not imply the usual Markov conclusion for the zeros. For deformations of the
first generator, \eqref{eq:nik-tail} is necessary and sufficient for the Markov
conclusion to hold for every common non-decreasing admissible representing
function, while
Theorems~\ref{thm:minimal-support}, \ref{thm:nik-christoffel},
and~\ref{thm:resolvent} cover the stated positive cases. More generally, at a
fixed parameter value, Theorem~\ref{thm:dual} gives the complete test for the
corresponding independent or common class of non-decreasing representing
functions: a negative tail is realised by the step deformation above.
The Angelesco theorem thus gives a general extension of Markov's principle,
while the dual criterion characterizes its validity for these deformation
classes in the normal type II setting with simple real zeros.

\section*{Acknowledgements}
This work was supported by the Portuguese Foundation for Science and Technology (FCT) under the project UID/00324/2025, Centre for Mathematics of the University of Coimbra (\href{https://doi.org/10.54499/UID/00324/2025}{doi:10.54499/UID/00324/2025}).  The first author acknowledges financial support from FCT under the project 2022.00143.\allowbreak CEECIND/\allowbreak CP1714/\allowbreak CT0002 (\href{https://doi.org/10.54499/2022.00143.CEECIND/CP1714/CT0002}{doi:10.54499/2022.00143.CEECIND/CP1714/CT0002}).  The second author acknowledges support from FCT through the grant UI.BD.154694. 2023.
\bibliographystyle{plain}
\bibliography{bib_revised}

\begin{thebibliography}{10}

\bibitem{A98}
A.~I. Aptekarev.
\newblock Multiple orthogonal polynomials.
\newblock {\em J. Comput. Appl. Math.}, 99(1--2):423--447, 1998.

\bibitem{AW85}
R.~Askey and J.~Wilson.
\newblock {\em Some basic hypergeometric orthogonal polynomials that generalize
  {J}acobi polynomials}, volume~54 of {\em Memoirs of the American Mathematical
  Society}.
\newblock American Mathematical Society, Providence, RI, 1985.

\bibitem{C22}
K.~Castillo.
\newblock Markov's theorem for weight functions on the unit circle.
\newblock {\em Constr. Approx.}, 55(2):605--627, 2022.

\bibitem{CZ22}
K.~Castillo and I.~Zaballa.
\newblock On a formula of {T}hompson and {M}c{E}nteggert for the adjugate
  matrix.
\newblock {\em Linear Algebra Appl.}, 634:37--56, 2022.

\bibitem{DaemsKuijlaars}
E.~Daems and A.~B.~J. Kuijlaars.
\newblock A {C}hristoffel--{D}arboux formula for multiple orthogonal
  polynomials.
\newblock {\em J. Approx. Theory}, 130(2):190--202, 2004.

\bibitem{D23}
A.~Doliwa.
\newblock Determinantal approach to multiple orthogonal polynomials and the
  corresponding integrable equations.
\newblock {\em Stud. Appl. Math.}, 153(2):e12726, 2024.

\bibitem{D17}
E.~J.~C. dos Santos.
\newblock Monotonicity of zeros of {J}acobi-{A}ngelesco polynomials.
\newblock {\em Proc. Amer. Math. Soc.}, 145(11):4741--4750, 2017.

\bibitem{FidalgoPerfect}
U.~Fidalgo~Prieto and G.~L\'opez~Lagomasino.
\newblock {N}ikishin systems are perfect.
\newblock {\em Constr. Approx.}, 34(3):297--356, 2011.

\bibitem{G71}
G.~Freud.
\newblock {\em Orthogonal polynomials}.
\newblock Pergamon Press, Oxford-New York, 1971.

\bibitem{I05}
M.~E.~H. Ismail.
\newblock {\em Classical and quantum orthogonal polynomials in one variable},
  volume~98 of {\em Encyclopedia of Mathematics and Its Applications}.
\newblock Cambridge University Press, Cambridge, 2005.

\bibitem{JST20}
C.~R. Johnson, R.~L. Smith, and M.~J. Tsatsomeros.
\newblock {\em Matrix positivity}, volume 221 of {\em Cambridge Tracts in
  Math.}
\newblock Cambridge University Press, Cambridge, 2020.

\bibitem{KY98}
A.~N. Kolmogorov and A.~P. Yushkevich, editors.
\newblock {\em Mathematics of the 19th century}, volume~3.
\newblock Birkh\"{a}user Verlag, Basel, 1998.
\newblock Function theory according to Chebyshev, ordinary differential
  equations, calculus of variations, theory of finite differences.

\bibitem{KozhanVaktnasChristoffel}
R.~Kozhan and M.~Vaktn\"as.
\newblock {C}hristoffel transform and multiple orthogonal polynomials.
\newblock {\em J. Comput. Appl. Math.}, 476:117121, 2026.

\bibitem{KozhanVaktnasZeros}
R.~Kozhan and M.~Vaktn\"as.
\newblock Zeros of multiple orthogonal polynomials: location and interlacing.
\newblock {\em Bull. Lond. Math. Soc.}, 58(1):e70281, 2026.

\bibitem{L22}
M.~Leurs.
\newblock {\em Jacobi-Angelesco multiple orthogonal polynomials and
  applications}.
\newblock PhD thesis, KU Leuven, Arenberg Doctoral School, 2022.

\bibitem{M95}
A.~P. Magnus.
\newblock Painlev{\'e}-type differential equations for the recurrence
  coefficients of semi-classical orthogonal polynomials.
\newblock {\em J. Comput. Appl. Math.}, 57(1--2):215--237, 1995.

\bibitem{M86}
A.~Markoff.
\newblock Sur les racines de certaines \'equations. {S}econde note.
\newblock {\em Math. Ann.}, 27:177--182, 1886.

\bibitem{MM24}
A.~Mart\'{\i}nez-Finkelshtein and R.~Morales.
\newblock Interlacing and monotonicity of zeros of {A}ngelesco-{J}acobi
  polynomials.
\newblock {\em Pure Appl. Funct. Anal.}, 9(5):1259--1279, 2024.

\bibitem{NS91}
E.~M. Nikishin and V.~N. Sorokin.
\newblock {\em Rational Approximations and Orthogonality}, volume~92 of {\em
  Translations of Mathematical Monographs}.
\newblock American Mathematical Society, Providence, RI, 1991.

\bibitem{SwiderskiVanAssche}
G.~{\'S}widerski and W.~Van Assche.
\newblock {C}hristoffel functions for multiple orthogonal polynomials.
\newblock {\em J. Approx. Theory}, 283:105820, 2022.

\bibitem{S75}
G.~Szeg\H{o}.
\newblock {\em Orthogonal polynomials}, volume~23 of {\em American Mathematical
  Society Colloquium Publications}.
\newblock American Mathematical Society, Providence, R. I., {F}ourth edition,
  1975.

\bibitem{Tao}
T.~Tao.
\newblock {\em Topics in random matrix theory}, volume 132 of {\em Grad. Stud.
  Math.}
\newblock American Mathematical Society, Providence, RI, 2012.

\end{thebibliography}

\end{document}